\documentclass[reqno,11pt]{amsart}
\usepackage[foot]{amsaddr}
\usepackage{amssymb,amsmath,amsthm,amsfonts}
\usepackage{mathrsfs,dsfont,comment,mathscinet,mathtools}
\usepackage{regexpatch}
\usepackage{graphicx,pgfplots}
\usepackage{cite}

\usepackage[utf8]{inputenc}
\usepackage[normalem]{ulem}
\usepackage[left=2cm,right=2cm,top=2.3cm,bottom=2.3cm]{geometry}

\usepackage{tkz-fct}

\usepackage{graphicx,tikz,color}
\usetikzlibrary{quotes,angles}

\usepackage[labelfont=bf]{caption}

\usepackage{enumerate}

\usepackage{ wasysym }

\usepackage{epstopdf}
\usepackage[colorlinks=true,  linkcolor=blue, 
            citecolor=blue, urlcolor=blue]{hyperref}

\allowdisplaybreaks
\numberwithin{equation}{section}
\theoremstyle{plain}
\newtheorem{theorem}{Theorem}[section]
\newtheorem{proposition}[theorem]{Proposition}
\newtheorem{lemma}[theorem]{Lemma}

\newtheorem{definition}[theorem]{Definition}

\theoremstyle{definition}

\newtheorem*{theorem*}{Theorem}

\newcommand{\Z}{\mathbb{Z}}
\newcommand{\N}{\mathbb{N}}
\newcommand{\R}{\mathbb{R}}
\newcommand{\RN}{\mathbb{R}^N}
\newcommand{\RNN}{\mathbb{R}^{N\times N}}
\newcommand{\SN}{\mathbb{S}^{N-1}}
\renewcommand{\H}{\mathscr{H}}
\renewcommand{\S}{\mathbb{S}}

\renewcommand{\d}{\mathrm{d}}
\newcommand{\restr}[1]{|_{#1}}

\usepackage{accents}
\renewcommand*{\dot}[1]{%
	\accentset{\mbox{\large\bfseries .}}{#1}}

\def\XXint#1#2#3{{\setbox0=\hbox{$#1{#2#3}{\int}$ }
\vcenter{\hbox{$#2#3$ }}\kern-.6\wd0}}

\makeatletter
\newcommand{\subalign}[1]{%
	\vcenter{%
		\Let@ \restore@math@cr \default@tag
		\baselineskip\fontdimen10 \scriptfont\tw@
		\advance\baselineskip\fontdimen12 \scriptfont\tw@
		\lineskip\thr@@\fontdimen8 \scriptfont\thr@@
		\lineskiplimit\lineskip
		\ialign{\hfil$\m@th\scriptstyle##$&$\m@th\scriptstyle{}##$\hfil\crcr
			#1\crcr
		}%
	}%
}
\makeatother

\makeatletter
\newcommand{\subsetcong}{\mathrel{\mathpalette\subseteq@cong\relax}}
\newcommand{\subsetsim}{\mathrel{\mathpalette\subset@sim\relax}}

\newcommand{\subseteq@cong}[2]{%
	\vbox{\offinterlineskip\m@th
		\ialign{\hfil$#1##$\hfil\cr
			\sim\cr\subseteq\cr
		}%
	}%
}

\newcommand{\subset@sim}[2]{%
	\vbox{\offinterlineskip\m@th
		\ialign{\hfil$#1##$\hfil\cr
			\sim\cr\subset\cr
		}%
	}%
}

\makeatother

\newcommand{\tr}{\mathrm{tr}}
\newcommand{\cof}{\mathrm{cof}}
\newcommand{\adj}{\mathrm{adj}}

\newcommand{\per}{\mathrm{Per}}
\renewcommand{\div}{\mathop{\mathrm{div}}}
\newcommand{\spt}{\mathop{\mathrm{spt}}}

\newcommand{\dist}{\mathrm{dist}}

\newcommand{\SO}{\mathrm{SO}(N)}

\newcommand{\SL}{\mathrm{SL}(N)}
\newcommand{\imt}{\mathrm{im}_{\mathrm{T}}}

\newcommand*\closure[1]{\overline{#1}}

\newcommand{\leb}{\mathscr{L}^N}
\newcommand{\haus}{\mathscr{H}^{N-1}}

\renewcommand{\O}{\mathrm{O}}

\newcommand{\id}{\mathrm{id}}
\newcommand{\Lip}{\mathrm{Lip}}

\definecolor{rr}{rgb}{0.7, 0.11, 0.11}
\definecolor{bb}{rgb}{0.2, 0.2, 0.6}

\makeatletter
\xpatchcmd{\proof}{\itshape}{\bfseries}{}{}
\makeatother

\title[Quasi-static cavity growth in incompressible elastomers]{On the quasi-static growth of cavities \\ in incompressible elastomers}
\author[M. Bresciani]{Marco Bresciani${}^*$}
\email{marco.bresciani@jku.at}
\author[G. Lazzaroni]{Giuliano Lazzaroni${}^\sharp$}
\email{giuliano.lazzaroni@unifi.it}
\address{${}^*$ Institute of Analysis, Johannes Kepler University Linz, Altenberger Strasse 69, 4040 Linz, Austria}
\address{${}^\sharp$ Department of Mathematics and Computer Science ``Ulisse Dini'', University of Florence, viale Morgagni 67/A, 50134 Firenze, Italy}

\date{\today}

\begin{document}

\setlength\parindent{0pt}

\vskip .2truecm
\begin{abstract}
For a variational model of incompressible elastomers exhibiting cavitation, we establish the existence of irreversible quasi-static evolutions driven by time-dependent boundary conditions. Our result extends previous work by Mora-Corral (2014) to incompressible materials under weaker regularity assumptions on the boundary datum. The proof employs the multiplicative splitting introduced by Francfort and Mielke (2006) to reduce the problem to  time-independent boundary conditions and then invokes abstract results on the existence of energetic solutions to rate-independent systems. Specifically, the construction of mutual recovery sequences requires particular attention to accommodate the incompressibility constraint:  to this end we employ  the   method by Moser (1965) to   exhibit  incompressible diffeomorphisms   interpolating between localized translations and the identity. 
\end{abstract}
\maketitle

\section{Introduction}\label{sec:intro}

\subsection{Motivation}
Cavitation is a failure mechanism typical of soft solids, which exhibit the   sudden  nucleation of holes in response to loading. Classically observed in rubbery elastomers subject to tensile stresses above a critical hydrostatic threshold \cite{gent.lindley57,gent.lindley59,gent.wang}, the phenomenon has also been reported in experiments on ductile metals~\cite{kassner.etal,petrinic1} and biological tissues \cite{dougan,zimberlin}. Its relevance to fracture is underscored by the long-established role of void growth and coalescence in the initiation of cracks \cite{rice.tracey}.

In rubber-like solids, resistance to volumetric deformation greatly exceeds resistance to shear, making incompressibility a natural idealization. In this setting, cavitation allows for macroscopic dilation through the formation of voids while preserving the local volume of the surrounding material \cite{lp.etal}.

Although incompressibility is a common simplifying assumption in the mechanical modeling of soft solids, the rigorous enforcement of this constraint often introduces  analytical subtleties of intrinsic mathematical interest.

\subsection{State of the art}
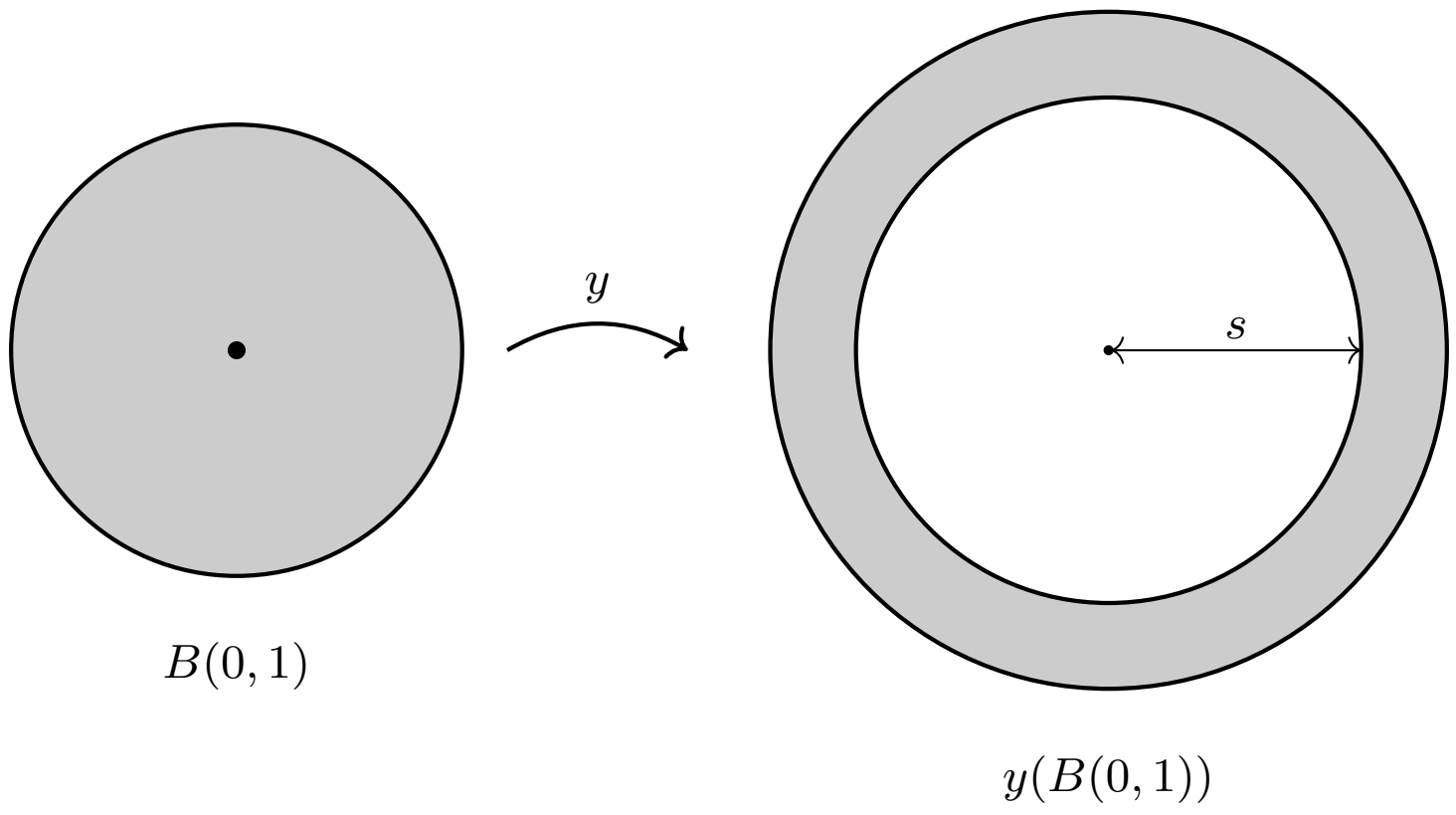
\begin{figure}
 	\centering
 	\begin{tikzpicture}[scale=1.5]
 		\pgfmathsetmacro{\a}{1/2}
 		
 		\begin{scope}[xshift=-60]
 			\draw[fill=black!20,thick] (0,0) circle (1);
 			\draw[fill=black] (0,0) circle (1pt);
 			\draw[->,thick] (1.2,0) to[out=30,in=150] node[above,midway] {\footnotesize ${y}$}  (2,0);
 			\node at (0,-1.4) {\footnotesize $B(0,1)$};
 		\end{scope}

 		\begin{scope}[xshift=50]
 			\draw[fill=black!20,thick] (0,0) circle (1.5);
 			\draw[fill=white,thick] (0,0) circle (1.12);
 			\draw[fill=black] (0,0) circle (.5pt);
 			\draw[<->] (.01,0) -- (1.12,0)  node[midway,above,inner sep=2pt] {\footnotesize $s$};
 			\node at (0,-1.9) {\footnotesize ${y}(B(0,1))$};
 		\end{scope}
 	\end{tikzpicture}	
 	\caption{An incompressible radial deformation}
 	\label{fig:inc}
 \end{figure}

The first analytical treatment of cavitation in solids is due to Ball~\cite{ball1982}. Restricting to radially symmetric deformations, he showed the existence of equilibrium solutions to the equations of nonlinear elastostatics in which a cavity forms at the center. For incompressible materials, any such deformation necessarily takes the form
\begin{equation*}
	y(x)
	\coloneqq\sqrt[N]{|x|^N+s^N}\frac{x}{|x|} \quad \text{for all $x\in B(0,1)\setminus \{0\}$,}
\end{equation*}
where $s>0$ denotes the radius of the cavity nucleated at the origin; cf.\ Figure~\ref{fig:inc}. Ball’s analysis interprets cavitation as an elastic instability, characterized through bifurcation of the associated boundary-value problem. This perspective has since been extended in many directions; see, e.g., the review in~\cite{horgan.polignone}.

A different modeling approach was taken in the seminal paper of Müller and Spector \cite{mueller.spector}. Without any symmetry assumption, they established the existence of equilibria given by solutions to free-discontinuity problems in which the location of flaw points is not prescribed. Beyond the usual elastic contribution, the energy comprises an additional term proportional to the perimeter of the deformed configuration. This approach aligns with the modern mechanical interpretation of cavitation as an irreversible failure phenomenon analogous to fracture, rather than an elastic instability \cite{bcllp,kumar.lopezpamies,lrclp}. Within the same functional framework, Sivaloganathan and Spector \cite{sivaloganathan.spector} proposed an alternative well-posed model with prescribed flaw points according to the purely elastic interpretation, subsequently studied via core-radius approximation in \cite{henao,sivaloganathan.spector.tilakraj}. The analysis of Müller and Spector was then specialized to neo-Hookean materials by Conti and De Lellis \cite{conti.delellis}, who highlighted a substantial obstacle for generalizing the existence result of \cite{mueller.spector} to their setting; we refer, e.g., to \cite{bhmcr.harmonicdipoles,bhmcr.axi,bhmcr.relataxion} for recent developments in this direction.

Building upon the analysis of \cite{conti.delellis}, the model of \cite{mueller.spector} was refined by Henao and Mora-Corral \cite{henao.moracorral.lusin}, who proposed a surface term proportional to the perimeter of the cavities created by the deformations.   Comparing with \cite{mueller.spector}, 
this choice is more natural, as it does not penalize stretching of the boundary. The resulting surface-energy functional corresponds to a particular instance of a more general functional on special functions of bounded variation that measures the new surface created by deformations \cite{henao.moracorral.fracture} and plays a key role in the weak continuity of Jacobian determinants \cite{henao.moracorral.invertibility}. An anisotropic variant of the model by Henao and Mora-Corral was investigated in \cite{bresciani}, while  core-radius and phase-field approximations have been studied in~\cite{bresciani.friedrich.core} and \cite{henao.moracorral.xu}, respectively.

An evolutionary version of the model in \cite{henao.moracorral.lusin} was advanced by Mora-Corral \cite{moracorral}, who established the existence of irreversible quasi-static evolutions. Irreversibility is incorporated through an additional energetic term that penalizes the nucleation of cavities, and the resulting formulation bears strong similarities with variational models for brittle fracture à la Griffith \cite{francfort.marigo}. The proof strategy follows the by-now-standard time-discretization scheme \cite{mielke.roubicek}; see, e.g.,  \cite{bresciani.friedrich,bresciani.stroffolini,dalmaso.lazzaroni,francfort.marigo,francfort.mielke}. This evolutionary model of Mora-Corral was subsequently extended to a two-dimensional framework allowing for both cavitation and fracture in~\cite{bresciani.friedrich}.

\subsection{Present contribution}
In this work, we specialize the analysis of \cite{moracorral} to incompressible materials and establish the existence of irreversible quasi-static evolutions driven by time-dependent boundary conditions. In doing so, we also weaken the regularity requirement on the boundary datum in line with the assumptions of \cite{lazzaroni} and, later, \cite{bresciani.stroffolini}.

To avoid repeating well-known arguments, we do not retrace all steps in \cite{moracorral}, but instead invoke the abstract existence theory for energetic solutions to rate-independent systems \cite{mielke.roubicek}. Since the incompressibility constraint does not affect lower-semicontinuity methods, our main effort is concentrated on proving the existence of mutual recovery sequences \cite[Eq.~(2.1.38)]{mielke.roubicek}. Alike in \cite{moracorral}, this requires establishing a cavity transfer result, given here in Proposition~\ref{prop:transfer}. Its proof hinges on the composition of deformations with suitable incompressible diffeomorphisms that relocate points, whose existence is established in Lemma~\ref{lem:incompressible}. The construction, which lies at the core of our analysis, combines the point-moving diffeomorphisms of \cite{moracorral} with a result on the prescribed Jacobian equation, also known as the pull-back equation~\cite{csato.dacorogna.kneuss}. Specifically, we employ the Moser flow method \cite{moser} to construct an incompressible diffeomorphism interpolating between localized translations and the identity. With this result at hand, the proof of Proposition~\ref{prop:transfer} proceeds largely as in \cite{moracorral}.

As in \cite{lazzaroni}, the boundary datum is assumed to be a time-by-time bi-Lipschitz homeomorphism with $W^{1,\infty}$-regularity in time. Although \cite{moracorral} anticipates the possibility of working under this assumption, doing so requires some technical adjustments. Above all stands the revision of the first-variation formula for the perimeter under bi-Lipschitz vector fields. Nonetheless, this regularity assumption is robust enough to allow one to carry it out  by means of standard tools in geometric measure theory 
\cite{ambrosio.fusco.pallara,maggi}.

As in \cite{moracorral}, our analysis does not cover the neo-Hookean model, since we must assume superlinear growth with respect to the cofactor for the elastic-energy density; see assumption~(W3). This limitation is not merely technical: it reflects the lack of compactness for finite-energy deformations first highlighted by Conti and De Lellis \cite{conti.delellis}. While a counterexample analogous to the one in \cite[Section~6]{conti.delellis} for incompressible materials is conceivable in principle, none has appeared in the literature to the best of our knowledge.

The present study originates from an attempt to extend the existence result in \cite{bresciani.friedrich} to incompressible elastomers. However, as noted in \cite[Footnote~2]{francfort.giacomini.lopezpamies}, an incompressible variant of the crack transfer poses a substantial mathematical obstacle that appears to lie beyond current technical reach.
\subsection{Structure of the paper}
 The paper comprises three sections and two appendices.  Section~\ref{sec:main}   recalls the relevant background,  lists all assumptions, and states our main result,  Theorem \ref{thm:main}. Section~\ref{sec:proof} is devoted to its proof. Appendix \ref{sec:aa} revisits the first-variation formula for the perimeter, and Appendix~\ref{sec:ab} summarizes the relevant results on the prescribed Jacobian equation.

\subsection{Notation and conventions}
For $r,s\in \R$, we set $r\wedge s \coloneqq \min\{r,s\}$ and $r\vee s \coloneqq \max\{r,s\}$.
We write $\R^{d\times d}$ for the space of real square matrices of order $d\in \N$. Moreover,
$\mathrm{SL}(d)$ and $\mathrm{SO}(d)$ denote the special linear group and
the special orthogonal group, respectively.

Throughout the paper, $N\in\N$ denotes the space dimension, which satisfies $N\geq 2$. 
For $A,B\in\RNN$, we use the Frobenius inner product $A:B \coloneqq \tr(A^\top B)$ and the associated norm $
|A| \coloneqq \sqrt{A:A}$.  
Given $v,w\in\RN$, we define $v\otimes w \coloneqq vw^\top \in \RNN$,
where vectors are identified with column matrices.

A {domain} is an open and connected subset of $\RN$. We use standard
notation for the interior,  closure, and boundary of a set. Given $E,F\subset\RN$, we write
$E\subset \subset F$ whenever $\overline{E}\subset F$ and $\overline{E}$ is compact. The support of a function is denoted by `$\operatorname{spt}$'.

We denote by $\leb$ and $\mathscr{H}^\alpha$ the Lebesgue measure on $\RN$ and 
 the $\alpha$-dimensional Hausdorff measure on $\RN$, respectively. We will only consider $\alpha\in\{0,N-1\}$. We will also employ the Lebesgue measure $\mathscr{L}^1$ on $\R$. Any statement
involving `almost everywhere' or similar terminology without any specification is understood with
respect to Lebesgue measure, unless differently specified.

We employ standard notation for function spaces. In particular, for
$k\in\N\cup\{0,\infty\}$, $C^k$ denotes the space of $k$ times continuously
differentiable functions, while $L^p$ and $W^{1,p}$, with
$p\in[1,\infty]$, denote Lebesgue and Sobolev spaces, respectively. The space
of Lipschitz functions is denoted by `$\Lip$'. The subscripts `$\mathrm{b}$' and
`$\mathrm{c}$' indicate bounded and compactly supported functions, respectively.
Whenever appropriate, maps in $W^{1,\infty}$ are identified with their
continuous representatives.

If $\mathrm{O}\subset\RN$ is open, we denote by $\mathrm{Diff}^k(\overline{\mathrm{O}};\overline{\mathrm{O}})$ 
the class of diffeomorphisms $\varphi\colon \mathrm{O}\to\mathrm{O}$ such
that both $\varphi$ and $\varphi^{-1}$ belong to
$C^k(\overline{\mathrm{O}};\RN)$. If $\Xi\subset\RN$ is compact, we write $W^{1,\infty}(\Xi;\Xi)$
in place of $W^{1,\infty}(\Xi^{\circ};\Xi)$.

Traces of Sobolev maps are denoted by restriction. The symbols $D$ and
`$\operatorname{div}$' denote the weak gradient and the weak divergence,
respectively. The same notation is used for weak tangential gradients and
divergences; in this case, the relevant rectifiable set is indicated by a
superscript. For these notions, we refer, e.g., to \cite[Appendix]{bresciani.friedrich.core}. The prime notation is reserved for weak
derivatives of functions defined on the real line.

We use the standard terminology of sets of finite perimeter as, e.g., in \cite{ambrosio.fusco.pallara, maggi}. For a measurable set
$E\subset\RN$, we write $\per(E)$ for its perimeter and $\partial^*E$ for its
reduced boundary.

 Given $T>0$ and a Banach space $X$,  we use the notation
$C^k([0,T];X)$ and $\Lip([0,T];X)$ for spaces of continuously differentiable
and Lipschitz maps, respectively, $L^p(0,T;X)$ for Bochner spaces, and
$W^{1,p}(0,T;X)$ for Bochner--Sobolev spaces. The subscripts `$\mathrm{b}$' and
`$\mathrm{c}$' retain the same meaning in this setting. As above, elements of
$W^{1,\infty}(0,T;X)$ are identified with their continuous representatives.
The time derivative of a Bochner--Sobolev map is denoted by a superposed dot.
For a map defined on $[0,T]$ with values in any abstract space, its value at time $t$ is  indicated by
a subscript $t$.

Finally, $C>0$ denotes a generic constant, whose value may change from line to
line. In some instances, its dependence on relevant quantities is highlighted by means of parentheses.

\section{Setting and main result} \label{sec:main}
In this section, we first describe the setting of the problem. Then, we state our main result with all the assumptions.  

\subsection{Variational modeling of cavitation}
In order to present our  findings, we briefly review the variational theory of cavitation in nonlinear elasticity under critical integrability. For more details, we refer to  \cite{conti.delellis,henao.moracorral.lusin,moracorral,mueller.spector}. 

For the rest of the paper, we fix 
\begin{equation*}
	\Omega\subset \RN \ \text{bounded Lipschitz domain}
\end{equation*}
representing the reference configuration of an elastomer. We enforce a physical confinement by fixing 
\begin{equation*}
	\Xi\subset \RN \  \text{closure of a bounded Lipschitz domain.}   
\end{equation*}
Deformations are modeled as maps in $W^{1,N-1}(\Omega;\Xi)$, i.e., $y\in W^{1,N-1}(\Omega;\RN)$ with $y(x)\in \Xi$ for almost all $x\in \Omega$. For any such map, a notion of degree is available.

\begin{definition}[Topological degree]\label{def:deg}
Let $y\in W^{1,N-1}(\Omega;\Xi)$ and $U\subset \subset \Omega$ domain of class $C^1$ with $y\restr{\partial U}\in W^{1,N-1}(\partial U;\RN)$. The {topological degree} of $y$ with respect to $U$ is the unique map 	$\deg(y,U,\cdot)\in\mathrm{BV}(\RN;\Z)$ characterized by the identity 
\begin{equation*}
	\int_{\RN} \deg(y,U,\xi)\div \psi(\xi)\,\d\xi=\int_{\partial U}\psi(y(x))\cdot (\cof D^{\partial U}y(x))\nu_U(x)\,\d\haus(x)
\end{equation*}
for all $\psi\in C^1(\RN;\RN)$. 
\end{definition}
The existence of such a function is proved in \cite[Remark~3.3]{conti.delellis}; see also \cite[Section~2.5]{dolezalova.hencl.maly}.
Notably, as observed in \cite[Remark~3.2]{conti.delellis},  $\deg(y,U,\cdot)$ coincides with the Brouwer degree for $y\restr{\partial U}\in C^0(\partial U;\RN)$, which is known to be continuous.
However,  the continuity of the topological degree is not guaranteed for $y\restr{\partial U}\in W^{1,N-1}(\partial U;\RN)$.
 
When $y\in C^1(\closure{\Omega};\Xi)$ and it is (everywhere) injective, then the right-hand side of the formula above can be recast as
\[
\int_{\partial U}\psi(y(x))\cdot (\cof D^{\partial U}y(x))\nu_U(x)\d\haus(x)
= \int_{y(\partial U)}\psi(\xi)\cdot \nu_{y( U)} (\xi)\d\haus(x) ,
\]
thanks to area formula.

Next, we introduce the concept of topological image.
Loosely speaking, the topological image of $U$ under $y$ is the region enclosed by the closed curve $y(\partial U)$. 

\begin{definition}[Topological image]\label{def:topim}
Let  $y\in W^{1,N-1}(\Omega;\Xi)$ and $U\subset \subset \Omega$ be as in Definition~\ref{def:deg}. The {topological image} of $U$ under y is defined as 
\begin{equation*}
	\imt(y,U)\coloneqq \{ \xi\in\RN:\:\deg(y,U,\xi)\neq 0  \}^{(1)}.
\end{equation*} 
\end{definition}
Here, $V^{(1)}$ denotes the set of points of density one with respect to a measurable set $V\subset \RN$, namely, those $\xi \in \RN$ for which
\begin{equation*}
	\lim_{r\to 0^+} \frac{\leb(V\cap B(\xi,r))}{\leb(B(\xi,r))}=1.   
\end{equation*}
Owing to the possible lack of continuity of $\deg(y,U,\cdot)$, the set $\imt(y,U)$ is generally not open. 
Considering the set of density points (i.e., the measure-theoretic interior) in Definition~\ref{def:topim} is just technically convenient; cf. \cite[Lemma~3.8]{conti.delellis}.  

Having defined the topological image, we can now recall the topological condition (INV) first introduced in \cite{mueller.spector}.

\begin{definition}[Invertibility condition]\label{def:INV}
A map $y\in W^{1,N-1}(\Omega;\Xi)$ is termed to satisfy {condition (INV)} whenever, for all $a\in \Omega$ and for almost all $r\in (0,r_a)$ with $r_a\coloneqq \dist(a,\partial \Omega)$, the following conditions hold:
\begin{enumerate}[\rm (i)]
	\item $y\restr{\partial B(a,r)}\in W^{1,N-1}(\partial B(a,r);\Xi)$;
	\item $y(x)\in\imt(y,B(a,r))$ for almost all $x\in B(a,r)$;
	\item $y(x)\notin \imt(y,B(a,r))$ for almost all $x\in \Omega\setminus B(a,r)$.
\end{enumerate}
\end{definition}
As the name suggests, condition (INV) is related to invertibility. Namely, (INV) implies almost everywhere injectivity 
whenever $\det Dy>0$ almost everywhere \cite[Lemma~3.4]{conti.delellis}.

\begin{figure}
	\centering
	\begin{tikzpicture}[>=Stealth, scale=1.5]
		\tikzset{deformationstyle/.style={
				fill=black!20,draw=black,line width=0.8pt,line join=round}}
		
		\draw[deformationstyle] (0,0) rectangle (3,2.5);
		\fill (1.3,.5) circle (1.5pt);
		\fill (1.6,1.8) circle (1.5pt);
		\fill (1,2) circle (1.5pt);
		\node at (-.5,2.2) { $\Omega$};
		\node at (.7, .4) {$C_y$};
		\node at (1.9,1.8) {$a$};
		
		\draw[->,thick] (3.5,1.2) to[out=20,in=160]
		node[midway,above] {$y$} (5.2,1.2);

		\begin{scope}[shift={(5.7,0)},xscale=0.666666667,shift={(-5.7,0)}]
			\draw[deformationstyle]
			(5.7,0.1)
			to[out=105,in=255] (5.9,2.3)
			to[out=30,in=180] (8.8,2.8)
			to[out=0,in=180] (10.7,2.4)
			to[out=-100,in=80] (10.35,1.3)
			to[out=-100,in=90] (10.1,0.2)
			to[out=180,in=0] (8.9,-0.2)
			to[out=180,in=-5] cycle;

			\draw[xshift=225pt,yshift=12pt,fill=white,thick,scale=1.2,
			shift={(.25,0)},xscale=1.5,shift={(-.25,0)}]
			plot[smooth cycle] coordinates
			{(.1,.1) (.5,.3) (.3,.6) (0,.4)};
			\draw[xshift=180pt,yshift=2pt,fill=white,thick,scale=1.2,
			shift={(1.45,0)},xscale=1.5,shift={(-1.45,0)}]
			plot[smooth cycle] coordinates
			{(1.7,1.6) (1.7,1.9) (1.2,1.8) (1.2,1.5) (1.4,1.4)};
			\draw[xshift=150pt,yshift=5pt,fill=white,thick,scale=1.2,
			shift={(1.2,0)},xscale=1.5,shift={(-1.2,0)}]
			plot[smooth cycle] coordinates
			{(1,1.4) (1.1,1.3) (1.4,1.4) (1.4,1.6) (1.1,1.7) (1,1.5)};
		\end{scope}
		\node at (9.6,2.2) { ${y}(\Omega)$};
		\node at (8.25,2) {$E(y,a)$};
	\end{tikzpicture}
	\caption{The cavitation points $a\in C_y$ in the reference configuration $\Omega$ and  the cavities $E(y,a)$ in the deformed configuration $y(\Omega)$.}
	\label{fig:cav}
\end{figure}
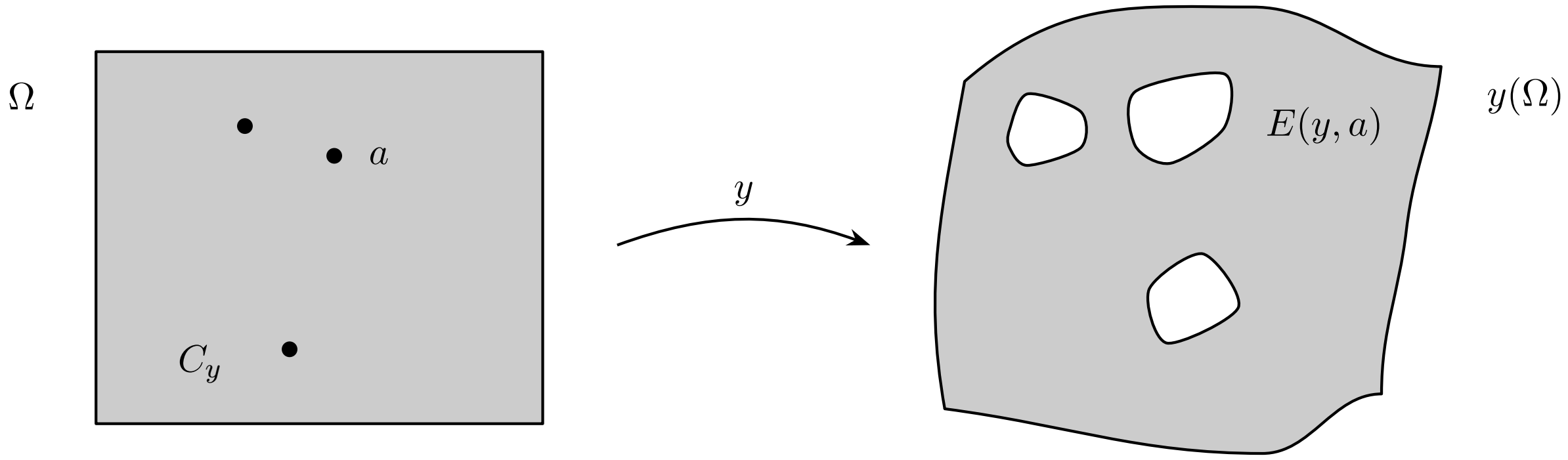

Next, we turn to the definition of cavities and cavitation points (cf.~Figure \ref{fig:cav}).

\begin{definition}[Cavity and cavitation points]
Let $y\in W^{1,N-1}(\Omega;\Xi)$ satisfy {\rm (INV)} and $\det D y>0$ almost everywhere. For any $a\in \Omega$, we define 
\begin{equation*}
	E(y,a)\coloneqq \bigcap \big \{ \imt(y,B(a,r)):\:r\in (0,r_a) \ \text{satisfying {\rm  (i)--(iii)} in Definition~\ref{def:INV}}  \big \}.
\end{equation*}
The set $C_y$ of {cavitation points} of $y$ is formed by those $a\in \Omega$ for which $\leb(E(y,a))>0$. For each $a\in C_y$, the set $E(y,a)$ is the {cavity} created by $y$ at $a$.
\end{definition} 
Note that $r\mapsto \imt(y,B(a,r))$ is increasing \cite[Lemma~3.12(ii)]{conti.delellis}. Also, more general sets than balls can be equivalently considered; cf.~\cite[Definition~2.18]{henao.moracorral.lusin}.

The following surface-energy functional has been introduced in \cite{henao.moracorral.invertibility} while studying the weak continuity of Jacobian determinants.

\begin{definition}[Surface-energy functional]\label{def:S}
Let $y\in W^{1,N-1}(\Omega;\Xi)$. Then, we define 
\begin{equation*}
	\mathcal{S}(y)\coloneqq \sup \left\{ \mathcal{S}(y,\Psi):\:\Psi\in C^\infty_{\rm c}(\Omega\times \RN;\RN),\ \|\Psi\|_{L^\infty({\Omega}\times \RN;\RN)}\leq 1  \right\},
\end{equation*}
where
\begin{equation*}
	\mathcal{S}(y,\Psi)\coloneqq \int_\Omega \Big(  D_x\Psi(x,y(x)):\cof Dy(x)+ \div_\xi \Psi(x,y(x))\det D y(x)  \Big)\,\d x
\end{equation*}
for any $\Psi \in C^\infty_{\rm c}(\Omega\times \RN;\RN)$ with variables $(x,\xi)\in \Omega\times \RN$.
\end{definition}
As shown in \cite{henao.moracorral.fracture}, the functional $\mathcal{S}$ correctly measures the area of the new surface created by deformations which, in the case of Sobolev mappings, correspond to the cavity boundaries. Precisely, the following holds \cite[Theorem~4.6]{henao.moracorral.lusin}: if $y\in W^{1,N-1}(\Omega;\Xi)$ with $\det D y>0$ almost everywhere satisfies condition (INV) and $\mathcal{S}(y)<+\infty$, then $C_y$ is countable and there holds
\begin{equation}
	\label{eq:S-cavities}
	\mathcal{S}(y)=\sum_{a\in C_y} \per(E(y,a)).
\end{equation}  
We remark that $E(y,a)\subset \Xi$ up to negligible sets in view of \cite[Lemma~2.5(iv)]{moracorral}.

\subsection{Main result}
Let $\Omega$ and $\Xi$ be as in the previous subsection. Additionally, assume that $\Omega\subset \Xi$. Incompressibility standing, deformations belong to the class
\begin{equation*}
	\mathcal{Y}\coloneqq \left\{ y\in W^{1,N-1}(\Omega;\Xi):\:\det D y=1 \ \text{a.e.,}\quad y \ \text{satisfies (INV),}\quad \mathcal{S}(y)<+\infty\right\},
\end{equation*} 
where we refer to Definition \ref{def:INV} and Definition~\ref{def:S}. We define the sets of flaw points as elements of
\begin{equation*}
	\mathcal{A}\coloneqq \left\{ A\subset \closure{\Omega}:\:\H^0(A)<+\infty   \right\},
\end{equation*}
and, in turn, the class of admissible configurations  reads
\begin{equation*}
	\mathcal{C}\coloneqq \left\{(y,A)\in\mathcal{Y}\times \mathcal{A}:\:C_y\subset A \right\}.
\end{equation*} 
To impose boundary conditions, we fix
\begin{equation*}
	\Gamma\subset \partial \Omega \ \text{rectifiable set with} \ \haus(\Gamma)>0.
\end{equation*}
For a time horizon $T>0$, we consider a boundary datum $w\colon [0,T]\to W^{1,\infty}(\Xi;\Xi)$  such that 
\begin{equation}
	\label{eq:w-datum}
	w_t \ \text{is invertible for all} \ t\in [0,T], \qquad w,w^{-1}\in W^{1,\infty}(0,T;W^{1,\infty}(\Xi;\Xi)); 
\end{equation}
cf.~\cite{lazzaroni}. Here and henceforth, we define $w^{-1}\colon t\mapsto  w_t^{-1}$. Additionally, we assume
\begin{equation*}
	\det D w_t=1 \ \text{a.e. in}\ \Omega \quad \text{for all} \ t\in [0,T].
\end{equation*}
For each $t\in [0,T]$,  the class of admissible configuration at time $t$ is given by
\begin{equation*}
	\mathcal{C}_t\coloneqq \left\{ (y,A)\in\mathcal{C}:\:y=w_t \ \text{on}\ \Gamma \right\},
\end{equation*}
where the equality is understood in the sense of traces. Note that each $\mathcal{C}_t$ is nonempty, since $(w_t,\emptyset)\in\mathcal{C}_t$.
We define the energy functional 
\begin{equation*}
	\mathcal{F}(y,A)\coloneqq \mathcal{G}(y)+\mathcal{H}(y)+\mathcal{I}(A),
\end{equation*}
where
\begin{equation*}
	\mathcal{G}(y)\coloneqq \int_\Omega W(Dy)\,\d x, \qquad \mathcal{H}(y)\coloneqq \sum_{a\in C_y} \kappa(a) \per \left(E(y,a) \right), \qquad \mathcal{I}(A)\coloneqq \sum_{a\in A} \lambda(a).
\end{equation*}
The functional $\mathcal{G}$ accounts for the elastic energy, while  $\mathcal{H}$ and $\mathcal{I}$ penalize the expansion of cavities and their nucleation, respectively. 

Here, $W\colon \RNN \to [0,+\infty]$ is the stored-energy density  which has the following properties: 
\makeatletter
\begin{enumerate}
	\renewcommand{\labelenumi}{(W\arabic{enumi})} 
	\renewcommand{\theenumi}{(W\arabic{enumi})}  
	\renewcommand{\p@enumi}{}                     
	
	\item\label{it:W-inc} \emph{Incompressibility:} $W(F)<+\infty$ if and only if $F\in \SL$;
	\item\label{it:W-fi} \emph{Frame indifference:} $W(QF)=W(F)$ for all $Q\in\SO$ and $F\in \SL$;
	\item\label{it:W-growth}\emph{Growth:} There exist $c_W>0$ and $g_W\colon [0,+\infty)\to [0,+\infty)$ Borel-measurable and superlinear at infinity (i.e., $\lim_{s\to +\infty}g_W(s)/s=+\infty$)
	such that
	\begin{equation*}
		W(F)\geq c_W |F|^{N-1} + g_W(|\cof F|)\quad \text{for all}\ F\in \SL;
	\end{equation*}
	\item\label{it:W-polyconvex} \emph{Polyconvexity:} There exists a convex function $\widehat{W}\colon \prod_{i=1}^{N-1} \mathrm{SL}\left(\binom{N}{i}\right)\to [0,+\infty)$ such that
	\begin{equation*}
		W(F)=\widehat{W}(\mathrm{M}(F))\quad \text{for all}\ F\in \SL,
	\end{equation*}
	   where $\mathrm{M}(F)$ collects all minors of $F$ up to order $N-1$ (e.g., $\mathrm{M}(F)=F$ when $N=2$, and $\mathrm{M}(F)=(F,\adj F)$ when $N=3$);
	\item\label{it:W-regular} \emph{Regularity:} $W$ is of class $C^1$ on an open neighborhood of $\SL$ in $\RNN$;
	\item\label{it:W-mandel} \emph{Control of the Mandel stress tensor:} There exists $\alpha_W> 0$ and $\beta_W\geq 0$ such that
	\begin{equation*}
		|F^\top DW(F)|\leq \alpha_W \left(W(F)+\beta_W \right) \quad \text{for all}\ F\in \SL;
	\end{equation*}
	\item\label{it:W-kirchhoff} \emph{Continuity of the Kirchhoff stress tensor:} There exists $\gamma_W\geq 0$ and, for each $\varepsilon>0$, there exists a $\delta_W=\delta_W(\varepsilon)>0$ such that
	\begin{equation*}
		|DW(GF)(GF)^\top -DW(F)F^\top|\leq \varepsilon \left( W(F)+\gamma_W\right) \quad \text{for all}\ F,G\in \SL
	\end{equation*}
	provided $|G-I|<\delta_W$.
\end{enumerate}
\makeatother
Condition \ref{it:W-inc} models the incompressibility of the material. The other assumptions \ref{it:W-growth}--\ref{it:W-kirchhoff} are by now standard; see, e.g.,~\cite{ball.op,dalmaso.lazzaroni,francfort.mielke,lazzaroni}. As shown in \cite[Proposition~2.1]{lazzaroni}, \ref{it:W-fi} and \ref{it:W-mandel} imply the control of the Kirchhoff stress tensor, namely
\begin{equation}
	\label{eq:control-kirchhoff}
		|DW(F)F^\top|\leq \alpha_W (W(F)+\beta_W) \quad \text{for all}\ F\in \SL,
\end{equation}
	where $\alpha_W>0$ and $\beta_W\geq 0$ are the constants in \ref{it:W-mandel}.

The densities $\kappa,\lambda\colon \closure{\Omega}\to [0,+\infty]$, which model possible inhomogeneities,  must satisfy:
\makeatletter
\begin{enumerate}
	\renewcommand{\labelenumi}{($\kappa\lambda$\arabic{enumi})} 
	\renewcommand{\theenumi}{($\kappa\lambda$\arabic{enumi})}  
	\renewcommand{\p@enumi}{}                     
	
	\item\label{it:kl-bdd} \emph{Boundedness from below:} There exists $c_{\kappa\lambda}>0$ such that 
	\begin{equation*}
		\kappa(a)\geq c_{\kappa\lambda}, \qquad \lambda(a)\geq c_{\kappa\lambda}\qquad \text{for all}\ a\in \closure{\Omega};
	\end{equation*} 
	\item \label{it:kl-continuous} \emph{Continuity:} $\kappa$ and $\lambda$ are continuous on $\closure{\Omega}$.
\end{enumerate}
\makeatother
Conditions \ref{it:kl-bdd}--\ref{it:kl-continuous} have  already been assumed in \cite{moracorral}. 
 By \ref{it:kl-bdd} and \eqref{eq:S-cavities}, we have  
 \begin{equation*}
 	\label{eq:HS}
 	\mathcal{H}(y)\geq c_{\kappa\lambda}\mathcal{S}(y).
 \end{equation*}
  From a modeling perspective, \ref{it:kl-bdd} implies that the creation of a cavity has always a positive energetic cost. Also, under this condition, we have $\H^0(A)<+\infty$ whenever $\mathcal{F}(y,A)<+\infty$, i.e., states with finite energy have at most finitely many flaw points. 
  Note that $\kappa$ and $\lambda$ might attain the value $+\infty$, thus allowing for prescribed flaw points; cf~\cite{sivaloganathan.spector}. 

The main result of the paper establishes the existence of irreversible quasi-static evolutions as in  \cite[Theorem~4.2]{moracorral} in the case of an  incompressible elastomer under weaker regularity assumptions on the boundary datum. 

\begin{theorem}[Existence of irreversible quasi-static evolutions]
	\label{thm:main}
Let $(y^0,A^0)\in\mathcal{C}_0$ be a stable initial datum, i.e.,  
\begin{equation*}
	\mathcal{F}(y^0,A^0)\leq \mathcal{F}(\widetilde{y},\widetilde{A})\quad \text{for all} \ (\widetilde{y},\widetilde{A})\in\mathcal{C}_0 \ \text{with} \ \widetilde{A}\supset A^0.
\end{equation*}
Then, there exists a measurable {irreversible quasi-static evolution}, namely, a function $(y,A)\colon [0,T]\to \mathcal{C}$ with $(y_t,A_t)\in\mathcal{C}_t$ for all $t\in [0,T]$ satisfying the following conditions:
\begin{enumerate}[\rm (i)]
	\item \emph{Initial condition:} $(y_0,A_0)=(y^0,A^0)$;
	\item \emph{Irreversibility:} $A$ is increasing, i.e.,  $A_s\subset A_t$ for all $s,t\in [0,T]$ with $s<t$;
	\item \emph{Minimality:} For each $t\in [0,T]$,
	\begin{equation*}
		\mathcal{F}(y_t,A_t)\leq \mathcal{F}(\widetilde{y},\widetilde{A}) \quad \text{for all}\ (\widetilde{y},\widetilde{A})\in\mathcal{C}_t \ \text{with}\ \widetilde{A}\supset A_t;
	\end{equation*}
	\item \emph{Power balance:} The function $t\mapsto \mathcal{F}(y_t,A_t)$ belongs to $W^{1,\infty}(0,T)$ and, for almost all $t\in (0,T)$,
	\begin{equation*}
		\begin{split}
			\frac{\d}{\d t}\mathcal{F}(y_t,A_t)&= \int_\Omega DW(Dy_t):D(\dot{w}_t\circ w_t^{-1}\circ y_t)\,\d x\\
			&\quad + \sum_{a\in C_{y_t}} \kappa(a) \int_{\partial^* E(y_t,a)} \div{}^{\partial^* E(y_t,a)}(\dot{w}_t\circ w_t^{-1})\,\d\haus.
		\end{split}
	\end{equation*}
\end{enumerate}
Moreover, one has $A_t=A^0\cup \bigcup_{s\in [0,t]} C_{y_s}$ for all $t\in [0,T]$. 
\end{theorem}

In Theorem \ref{thm:main}, the measurability of  $(y,A)$ is meant with respect to the product topology on $\mathcal{Y}\times \mathcal{A}$ when $\mathcal{Y}$ and $\mathcal{A}$ are endowed with the weak topology of $W^{1,N-1}(\Omega;\RN)$ and the topology induced by the Hausdorff distance on compact subsets of $\closure{\Omega}$, respectively.

\section{Proof of Theorem \ref{thm:main}}\label{sec:proof}

This section is devoted to the proof of Theorem \ref{thm:main}.
Following \cite{moracorral}, we reduce to time-independent boundary condition by using a multiplicative splitting  \cite{francfort.mielke}; see also for similar approaches \cite{bresciani.stroffolini,dalmaso.lazzaroni,lazzaroni}. Then, instead of following all the steps in \cite{moracorral}, we resort to the abstract theory of rate-independent systems \cite{mielke.roubicek} to prove the existence of irreversible quasi-static evolutions.

\subsection{Reduction to time-independent boundary conditions}
\label{sec:ti-bs} Following the approach in \cite{francfort.mielke}, we reduce to a formulation with time-independent boundary conditions. Specifically, we consider time-dependent deformations of the form $t\mapsto w_t\circ z$ with $z\in\mathcal{Y}$ satisfying $z=\id$ on $\Gamma$. Letting 
\begin{equation*}
	\widehat{\mathcal{C}}\coloneqq \left\{ (z,A)\in\mathcal{C}:\:z=\id \ \text{on} \ \Gamma \right\},
\end{equation*}
 the map $(z,A)\mapsto (w_t\circ z,A)$ provides a homeomorphism between $\widehat{\mathcal{C}}$ and $\mathcal{C}_t$, where both classes are equipped with the subspace topology of $\mathcal{C}$. This can be proved through the same arguments in \cite[Proposition~5.1]{moracorral}, up to minor technical adjustments. Specifically, given
our weaker regularity assumption \eqref{eq:w-datum} on the boundary datum, showing that $\mathcal{S}(w_t\circ z)<+\infty$ requires  an approximation by mollification as done in  \cite[Proposition~5.4]{bresciani.stroffolini}. The incompressibility of $w_t\circ z$ is promptly checked via the chain rule recalling that $\det Dw_t=1$ almost everywhere. 

We define the time-dependent functionals
\begin{equation*}
	\widehat{\mathcal{F}}(t,z,A)\coloneqq \mathcal{F}(w_t\circ z,A), \quad \widehat{\mathcal{G}}(t,z)\coloneqq \mathcal{G}(w_t\circ z), \quad \widehat{\mathcal{H}}(t,z)\coloneqq \mathcal{H}(w_t\circ z),
\end{equation*}
so that $\widehat{\mathcal{F}}(t,z,A)=\widehat{\mathcal{G}}(t,z)+\widehat{\mathcal{H}}(t,z)+\mathcal{I}(A)$. The main result for the problem with time-independent boundary conditions reads as follows. 

\begin{theorem}[Auxiliary existence result]\label{thm:auxiliary}
Let $(z^0,A^0)\in\widehat{\mathcal{C}}$ satisfy
\begin{equation*}
	\label{eq:id}
	\widehat{\mathcal{F}}(0,z^0,A^0)\leq \widehat{\mathcal{F}}(0,\widetilde{z},\widetilde{A}) \quad \text{for all}\ (\widetilde{z},\widetilde{A})\in\widehat{\mathcal{C}} \ \text{with} \ \widetilde{A}\supset A^0.
\end{equation*}
 Then, there exists a measurable function $(z,A)\colon [0,T]\to \widehat{\mathcal{C}}$ satisfying the following conditions:
\begin{enumerate}[\rm (i)]
	\item $(z_0,A_0)=(z^0,A^0)$;
	\item The function $t\mapsto A_t$ is increasing, i.e.,  $A_s\subset A_t$ for all $s,t\in [0,T]$ with $s<t$;
	\item For each $t\in [0,T]$, 
	\begin{equation*}
		\widehat{\mathcal{F}}(t,z_t,A_t)\leq \widehat{\mathcal{F}}(t,\widetilde{z},\widetilde{A})\quad \text{for all}\ (\widetilde{z},\widetilde{A})\in\widehat{\mathcal{C}}\ \text{with}\ \widetilde{A}\supset A_t;
	\end{equation*}
	\item The function $t\mapsto \widehat{\mathcal{F}}(t,z_t,A_t)$ belongs to $W^{1,\infty}(0,T)$ and, for almost all $t\in [0,T]$,
	\begin{equation*}
		\begin{split}
			\frac{\d}{\d t}\widehat{\mathcal{F}}(t,z_t,A_t)&= \int_\Omega DW(D(w_t\circ z_t)):D(\dot{w}_t\circ z_t)\,\d x\\
			&\quad+ \sum_{a\in C_{z_t}} \kappa(a)  \int_{w_t(\partial^*E(z_t,a))} \div{}^{w_t(\partial^*E(z_t,a))}(\dot{w}_t\circ w_t^{-1})\,\d\haus. 
		\end{split}
	\end{equation*}
\end{enumerate}
Moreover, one has $A_t=A^0\cup\bigcup_{s\in [0,t]}C_{z_s}$ for all $t\in [0,T]$.
\end{theorem}

Theorem \ref{thm:main} follows immediately from Theorem~\ref{thm:auxiliary}: choosing $z^0\coloneqq w_0^{-1}\circ y^0$ and defining $y_t\coloneqq w_t\circ z_t$, all claims are verified by changing variables with the aid of \cite[Proposition~6.1]{moracorral}.

\subsection{Reformulation in the framework of rate-independent systems}\label{sec:ris}
Henceforth, we will employ the terminology in \cite{mielke.roubicek}.  Consider the {rate-independent system} $(\mathcal{Q},\mathcal{E},\mathcal{D})$  given by  the  state space $\mathcal{Q}\coloneqq \mathcal{Y}\times \mathcal{A}$, and the   energy functional $\mathcal{E}\colon [0,T]\times \mathcal{Q}\to [0,+\infty]$ and  dissipation distance $\mathcal{D}\colon \mathcal{A}\times \mathcal{A}\to [0,+\infty]$  defined as
\begin{align} \label{def:E}
	 \mathcal{E}(t,z,A)&\coloneqq \begin{cases}
	 	\widehat{\mathcal{G}}(t,z)+\widehat{\mathcal{H}}(t,z) & \text{if} \ (z,A)\in\widehat{\mathcal{C}},\\
	 	+\infty & \text{otherwise,}
	 \end{cases}\\ \label{def:D}
	 \mathcal{D}(A,\widetilde{A})&\coloneqq \begin{cases}
		\sum_{a\in \widetilde{A}\setminus A} \lambda(a) & \text{if}\ A\subset \widetilde{A}, \\ +\infty & \text{otherwise.}
	\end{cases}
\end{align}
With these definitions, we have that a function $(z,A)\colon [0,T]\to \widehat{\mathcal{C}}$ fulfills conditions (ii)--(iv) in Theorem~\ref{thm:auxiliary} if and only if the function $q\colon [0,T]\to \mathcal{Q}$ given by $q_t\coloneqq (z_t,A_t)$ is an {energetic solution} to  $(\mathcal{Q},\mathcal{E},\mathcal{D})$. That is,  the following two conditions hold:
\begin{enumerate}[\rm (i)]
	\item \emph{Global stability:} 
	\begin{equation*}
		\mathcal{E}(t,q_t)\leq \mathcal{E}(t,\widetilde{q})+\mathcal{D}(A_t,\widetilde{A})\quad \text{for all}\ t\in [0,T]\ \text{and}\  \widetilde{q}=(\widetilde{A},\widetilde{z})\in\mathcal{Q};
	\end{equation*}
	\item \emph{Energy-dissipation balance:}
	\begin{equation*}
		\mathcal{E}(t,q_t)+\mathrm{Var}_{\mathcal{D}}(A;[0,t])=\mathcal{E}(0,q_0)+\int_0^t \partial_t\mathcal{E}(s,q_s)\,\d s \quad \text{for all}\ t\in [0,T]. 
	\end{equation*}
\end{enumerate} 
In (ii),  the variation of $A$ with respect to $\mathcal{D}$ is defined as
\begin{equation*}
	\mathrm{Var}_{\mathcal{D}}(A;[0,t])=\sup \left \{ \sum_{i=1}^{m} \mathcal{D}(A_{t_{i-1}},A_{t_i}):\:0=t_0<t_1\dots<t_m=t\right \}.
\end{equation*}
whereas $\partial_t \mathcal{E}(s,q)$ indicates the derivative of the function $t\mapsto \mathcal{E}(t,q)$ for $t=s$.   
The equivalence is readily established (cf. \cite[Section~4.2.4.1]{mielke.roubicek}). The global stability of $q$ corresponds to (iii) in Theorem~\ref{thm:auxiliary} since $\widehat{\mathcal{F}}(t,z_t,A_t)=\mathcal{E}(t,q_t)+\mathcal{D}(\emptyset, A_t)$ and $\mathcal{D}(\emptyset,\widetilde{A})=\mathcal{D}(\emptyset,A_t)+\mathcal{D}(A_t,\widetilde{A})$. Noting that $\mathrm{Var}_{\mathcal{D}}(A;[0,t])=\mathcal{D}(A_0,A_t)$ whenever $A$ is increasing in $[0,t]$, the equivalence between (iii)--(iv) in Theorem~\ref{thm:auxiliary} and the energy-dissipation balance is immediate after checking that the time derivative of the energy has the right expression (cf.~Proposition~\ref{prop:td-energy}).

In view of this observation, proving Theorem~\ref{thm:auxiliary} reduces to showing the existence of an energetic solution to $(\mathcal{Q},\mathcal{E},\mathcal{D})$. To this end, we apply \cite[Theorem~2.1.6]{mielke.roubicek}. Precisely, in Section~\ref{sec:sa}, we verify that $(\mathcal{Q},\mathcal{E},\mathcal{D})$ fulfills the main assumptions for rate-independent systems \cite[Eqs.~(D1)--(D2) \& (E1)--(E2)]{mielke.roubicek}, together with a suitable estimate for approximating the time derivative of the energy by incremental ratios \cite[Eq.~(2.1.27)]{mielke.roubicek}. Then, in Section~\ref{sec:ct}, we  establish the transfer of cavities (cf.~Proposition~\ref{prop:transfer}) and, following \cite[Section~4.2.4.1]{mielke.roubicek},  we note that this yields the existence of {mutual recovery sequences} \cite[Eq.~(2.1.38)]{mielke.roubicek}. Consequently, \cite[Corollary~2.1.19(a)]{mielke.roubicek} implies that $(\mathcal{Q},\mathcal{E},\mathcal{D})$ satisfies the two {compatibility conditions} \cite[Eqs.~(C1)--(C2)]{mielke.roubicek}. The existence of energetic solutions therefore follows from \cite[Theorem~2.1.6]{mielke.roubicek}. 

Eventually, one has to prove the characterization of flaw points in terms of the cavitation points of deformations along the evolution, which we do in Section~\ref{sec:flaw}.

\subsection{Structural assumptions}\label{sec:sa}
First, we check that $(\mathcal{Q},\mathcal{E},\mathcal{D})$ satisfies the main assumptions of rate-independent systems.

 The spaces $\mathcal{Y}$ and $\mathcal{A}$ are endowed with the weak topology of $W^{1,N-1}(\Omega;\RN)$ and the topology induced by the Hausdorff distance on compact subsets of $\closure{\Omega}$, respectively. Then,  
  $\mathcal{Q}$ is endowed with the product topology and, in turn, its compact subsets are both metrizable and separable.

The dissipation distance $\mathcal{D}$ in \eqref{def:D} is easily seen to satisfy the triangle inequality and $\mathcal{D}(A,\widetilde{A})=0$ if and only if $A=\widetilde{A}$ owing to \ref{it:kl-bdd}. The lower semicontinuity of $\mathcal{D}$ is checked below.

\begin{lemma}[Lower semicontinuity of the dissipation distance]
Let $(A_n)_{n\in\N}$ and $(\widetilde{A}_n)_{n\in\N}$ be  sequences in $\mathcal{A}$, and let $A,\widetilde{A}\in\mathcal{A}$ such that $A_n\to A$ and $\widetilde{A}_n\to \widetilde{A}$ in $\mathcal{A}$. Then,
\begin{equation} \label{260924}
	\mathcal{D}(A,\widetilde{A})\leq \liminf_{n\to \infty} \mathcal{D}(A_n,\widetilde{A}_n).
\end{equation}	
\end{lemma}
\begin{proof}
	Without loss of generality, we can assume that $A_n\subset \widetilde{A}_n$ for all $n\in\N$ and, by \ref{it:kl-bdd}, also 
	\begin{equation*}
		0<\sup_{n\in\N}\H^0(\widetilde{A}_n)<+\infty ,	
	\end{equation*}
otherwise the right-hand side of \eqref{260924} is infinite.
	As in \cite[Lemma~3.1]{bresciani.friedrich}, up to subsequences, we find $m\in\N\cup\{0\}$ and $M\in \N$ such that
$\H^0(A_n)=m$ and $\H^0(\widetilde{A}_n)=M$ for all $n$. It turns out that $\H^0(A)\le m$ and $\H^0(\widetilde{A})\le M$.
We may introduce enumerations of the elements of $A_n$,  $\widetilde{A}_n$, $A$, $\widetilde{A}$,
	\begin{equation*}
		A_n=\{a_{1,n},\dots,a_{m,n}\}, \quad \widetilde{A}_n=\{\widetilde{a}_{1,n},\dots,\widetilde{a}_{M,n}  \}, 
 \quad \text{without  repetitions,}
	\end{equation*}
\[ A=\{a_{1},\dots,a_{m}\}, \quad \widetilde{A}=\{\widetilde{a}_{1},\dots,\widetilde{a}_{M}\}, \quad \text{with possible repetitions,}
\]
in such a way that 
$a_{i,n}\to a_i$ and  $\widetilde{a}_{j,n}\to \widetilde{a}_{j}$, as $n\to \infty$, for all $i=1,\dots,m$ and $j=1,\dots,M$. Recalling $A_n\subset \widetilde{A}_n$, up to reordering we can assume $\widetilde{a}_{i,n}=a_{i,n}$ for $i=1,\dots,m$.  Thus, $A\subset \widetilde{A}$ and, by \ref{it:kl-continuous}, 
	\begin{equation*}
		\mathcal{D}(A,\widetilde{A})=\sum_{a\in \widetilde{A}\setminus  A} \lambda(a)\leq \sum_{i=m+1}^{M}\lambda(\widetilde{a}_i)= \lim_{n\to \infty} \sum_{i=m+1}^{M}\lambda(\widetilde{a}_{i,n})=\lim_{n\to \infty} \mathcal{D}(A_n,\widetilde{A}_n). \qedhere
	\end{equation*}
\end{proof}

We next turn to the energy functional $\mathcal E$ from \eqref{def:E}. By \ref{it:W-inc}, we have $\mathcal{E}(t,q)<+\infty$ if and only if $\mathcal{E}(0,q)<+\infty$. Moreover, for every $q\in\mathcal{Q}$, the functional $q\mapsto\mathcal{E}(t,q)$ is coercive and lower semicontinuous (equivalently, its sublevel sets are compact). These two properties can be proved as in \cite[Proposition~6.4 \& Proposition~6.5]{moracorral} by making use of \ref{it:W-growth}--\ref{it:W-polyconvex} and \ref{it:kl-bdd}--\ref{it:kl-continuous}.

The following result provides an energetic control of the power. Compared with \cite{moracorral}, our analysis is more delicate owing to the weaker regularity assumption \eqref{eq:w-datum} on the boundary datum. While the elastic energy can be treated as in \cite{lazzaroni}, the cavitation energy requires revisiting the first-variation formula for the perimeter functional. This task is deferred to Appendix \ref{sec:aa}.

\begin{proposition}[Energetic control of the power]\label{prop:td-energy}
	For each $(z,A)\in\widehat{\mathcal{C}}$, the function $t\mapsto\mathcal{E}(t,z,A)$ belongs to $W^{1,\infty}(0,T)$. Moreover, there exists $S_w\subset (0,T)$ with $\mathscr{L}^1(S_w)=0$ for which the following holds:
	\begin{enumerate}[\rm (i)]
		\item For all $t\in (0,T)\setminus S_w$, there exists 
		\begin{align}
			\label{eq:E-td}
			\begin{split}
				\partial_t{\mathcal{E}}(t,z,A) &=\int_\Omega DW(D(w_t\circ z)):D(\dot{w}_t\circ z)\,\d x\\
				&\quad +\sum_{a\in C_z} \kappa(a) \int_{w_t(\partial^*E(z,a))}\div{}^{w_t(\partial^*E(z,a))}(\dot{w}_t\circ w_t^{-1})\,\d\haus 
			\end{split}
		\end{align}
		with
		\begin{equation}
			\label{eq:E-td-est}
			|\partial_t \mathcal{E}(t,z,A)|\leq C(N,\Omega,W,w) \left  (\mathcal{E}(t,z,A)+1\right).
		\end{equation}
		\item For every $\varepsilon,c>0$ and for all $t\in (0,T)\setminus S_w$, there exists $\delta=\delta(\varepsilon,c,t)>0$ such that
		 \begin{equation}
		 	\label{eq:E-td-fo}
		 	\left | \frac{\mathcal{E}(t+h,z,A)-\mathcal{E}(t,z,A)}{h}-\partial_t\mathcal{E}(t,z,A) \right |\leq \varepsilon \quad \text{for all}\ h\in (-\delta,\delta)
		 \end{equation}
		 provided $\mathcal{E}(0,z,A) \leq c$.
	\end{enumerate}	
\end{proposition}
\begin{proof}
	For brevity, let $f(t)\coloneqq \widehat{\mathcal{F}}(t,z)$ and $g(t)\coloneqq \widehat{\mathcal{G}}(t,z)$, so that $\mathcal{E}(t,z,A)=f(t)+g(t)$. Let $S_w\subset (0,T)$ with $\mathscr{L}^1(S_w)$ be  such that $w$ is $W^{1,\infty}(\Xi;\Xi)$-differentiable at all $t\in (0,T)\setminus S_w$. The proof is divided into three steps.
	
	\emph{Step 1: Elastic energy.} 
	 Define $\widehat{W}\colon [0,T]\times \Xi \times \RNN\to [0,+\infty]$ by setting $\widehat{W}(t,\xi,F)\coloneqq W(D w_t(\xi)F)$. By the chain rule, 
	\begin{equation*}
		f(t)=\int_\Omega \widehat{W}(t,z(x),Dz(x))\,\d x \quad \text{for all}\ t\in [0,T].
	\end{equation*}
	 As in \cite[Remark~1.10]{lazzaroni}, we find that there exists $R_w\subset \Xi$ with $\leb(R_w)=0$ such that $w_t$ is differentiable at all $\xi\in \Xi\setminus R_w$ for all $t\in [0,T]$. Then, by \ref{it:W-regular} and \eqref{eq:w-datum}, we compute
	\begin{equation*}
		\partial_t \widehat{W}(t,\xi,F)=DW(Dw_t(\xi)F):(D\dot{w}_t(\xi)F ) \quad \text{for all}\ t\in (0,T)\setminus S_w,\ \xi\in \Xi\setminus R_w,\ \text{and}\ F\in\SL.
	\end{equation*} 
	With the same arguments in \cite[p.~181 \& Proposition~2.4]{lazzaroni}, which make use of \ref{it:W-mandel}--\ref{it:W-kirchhoff}, we thus obtain
	\begin{equation}
		\label{eq:f-dot}
		{f}'(t)=\int_\Omega \partial_t \widehat{W}(t,z(x),Dz(x))\,\d x\quad \text{for all}\ t\in (0,T)\setminus S_w.
	\end{equation}
	Writing $K(F)\coloneqq DW(F)F^\top$ for all $F\in \SL$, we find
	\begin{equation}
		\label{eq:W-overline-dot}
		\partial_t \widehat{W}(t,\xi,F)=K(Dw_t(\xi)F):D\dot{w}_t(\xi)(Dw_t(\xi))^{-1} \quad \text{for all}\ t\in (0,T)\setminus S_w.
	\end{equation}
	Thus, by \eqref{eq:w-datum} and \eqref{eq:control-kirchhoff},
	\begin{equation}
		\label{eq:f-dot-est}
		\begin{split}
			|{f}'(t)|&\leq \int_\Omega |K(Dw_t\circ zDz)|\,|Dw_t\circ z|\,|(Dw_t)^{-1}\circ z|\,\d x
			\leq C(w)\left\{ \alpha_W f(t) + \beta_W \leb(\Omega) \right\}.
		\end{split}
	\end{equation} 
	The proof in \cite{lazzaroni} also shows that, given $t\in (0,T)\setminus S_w$,
	\begin{equation*}
		\left | \frac{f(t+h)-f(t)}{h}-{f}'(t) \right |\leq \omega_{t}(h) \left(\alpha_W f(t)+\beta_W \leb(\Omega) \right) \quad \text{for all}\ h\in (t,T-t),  
	\end{equation*}
	where $\omega_t$ satisfies $\omega_t(h)\to 0$, as $h\to 0$.
	Gronwall's inequality and \eqref{eq:f-dot-est} yield
	\begin{equation*}
		f(t)\leq C(\Omega,W,w) (f(0)+1)\mathrm{e}^{C(W,w)t}\leq  C(\Omega,T,W,w) (f(0)+1),
	\end{equation*}
	so that
	\begin{equation}
		\label{eq:ffdot}
		\left | \frac{f(t+h)-f(t)}{h}-{f}'(t) \right |\leq \omega_{t}(h) C(\Omega,T,W,w)\left( f(0)+1 \right) \quad \text{for all}\ h\in (t,T-t).  
	\end{equation}
	
	\emph{Step 2: Cavitation energy.} We write $	g(t)=\sum_{a\in C_z} \kappa(a) P_{z,a}(t)$, where
	\begin{equation*}
		P_{z,a}(t)\coloneqq \per(w_t(E(z,a)))\quad \text{for all}\ t\in [0,T] \ \text{and}\ a\in C_z.	
	\end{equation*}
	By Theorem~\ref{thm:per}, we find that  $g\in W^{1,\infty}(0,T)$ with
	\begin{equation}
		\label{eq:g-dot}
		{g}'(t)=\sum_{a\in C_z} \kappa(a)\int_{w_t(\partial^*E(z,a))} \div{}^{w_t(\partial^*E(z,a))}(\dot{w}_t\circ w_t^{-1})\,\d\haus \quad \text{for all}\ t\in (0,T)\setminus S_w
	\end{equation}
	and
	\begin{equation}
		\label{eq:g-dot-est}
		|{g}'(t)|\leq C(N,w)\sum_{a\in C_z} \kappa(a) \per(E(z,a))\leq C(N,w) g(t).
	\end{equation}
	The last inequality is justified by observing that 
	\begin{equation}
		\label{eq:pp} \begin{split}
	 \per(E(z,a)) &  =\per(w_t^{-1}(w_t(E(z,a)))) \\ & \leq \| D w_t ^{-1}\|_{L^\infty(\Xi;\RNN)}^{N-1} \per(w_t(E(z,a)))\leq C(w) \per(w_t(E(z,a))) \end{split}
	\end{equation} 
	thanks to \eqref{eq:w-datum}.
	Additionally, Theorem~\ref{thm:per}  shows that there exists $\delta=\delta(w)$ and, for each $t\in (0,T)\setminus S_w$, a function $\omega_t\colon (-\delta,\delta)\to [0,+\infty)$ with $\omega_t(h)\to 0$, as $h\to 0$, such that, for all $h\in (-\delta,\delta)$,
	\begin{equation}
		\label{eq:ggdot}
		\begin{split}
			\left | \frac{g(t+h)-g(t)}{h}-{g}'(t)\right |&\leq C(N,w) \sum_{a\in C_z} \per(E(z,a)) \omega_t(h)\\
			&\leq C(N,w) \sum_{a\in C_z} \per(w_0(E(z,a))) \omega_t(h) .
		\end{split}
	\end{equation}
	For the latter inequality, we argue as in  \eqref{eq:pp} for $t=0$.
	
	\emph{Step 3: Conclusion.} Formula \eqref{eq:E-td} follows immediately from \eqref{eq:f-dot}--\eqref{eq:W-overline-dot} and \eqref{eq:g-dot}, while \eqref{eq:f-dot-est} and \eqref{eq:g-dot-est} yield \eqref{eq:E-td-est}. Eventually, the combination of \eqref{eq:ffdot} and \eqref{eq:ggdot} yields \eqref{eq:E-td-fo}.
\end{proof}

\subsection{Cavities transfer}\label{sec:ct}
In this subsection, we prove the cavity transfer for incompressible materials (cf.~\cite[Proposition~10.1]{moracorral}). Here, compared with \cite{moracorral}, the proof requires some extra care  to accommodate the incompressibility constraint. Also note that, in contrast with \cite[Proposition~10.1]{moracorral}, the coalescence of cavities is allowed in our assertion.

\begin{proposition}[Cavities transfer]\label{prop:transfer}
Let $(t_n)_{n\in\N}$ and $(A_n)_{n\in\N}$ be sequences in $[0,T]$ and $\mathcal{A}$, respectively, and let $t\in [0,T]$ and  $A\in \mathcal{A}$ such that $t_n\to t$ and  $A_n\to A$ in $\mathcal{A}$. Also,  let $(\widetilde{A},\widetilde{z})\in\widehat{\mathcal{C}}$ with $\widetilde{A}\supset A$ and $\widehat{\mathcal{F}}(t,\widetilde{z})<+\infty$.  Then, there exists a sequence $(\widetilde{A}_n,\widetilde{z}_n)_{n\in\N}$ in $\widehat{\mathcal{C}}$ with $\widetilde{A}_n\supset A_n$ for each $n\in\N$ such that
\begin{equation}
	 \label{eq:tc}
	\widetilde{z}_n\to \widetilde{z} \ \text{in}\ W^{1,N-1}(\Omega;\RN),  \qquad \widetilde{A}_n\to \widetilde{A} \ \text{in}\ \mathcal{A},
\end{equation}
and, additionally,
\begin{equation}
	\label{eq:tGH}
	\lim_{n\to \infty} \widehat{\mathcal{G}}(t_n,\widetilde{z}_n)=\widehat{\mathcal{G}}(t,\widetilde{z}), \qquad \lim_{n\to \infty} \widehat{\mathcal{H}}(t_n,\widetilde{z}_n)=\widehat{\mathcal{H}}(t,\widetilde{z}),
\end{equation}
\begin{equation}
	\label{eq:tI}
	\widetilde{A}_n\setminus A_n=\widetilde{A}\setminus A \ \text{for all} \ n\in\N.
\end{equation}
\end{proposition}

Before presenting the proof, we observe that the preceding result guarantees the existence of {mutual recovery sequences} for the rate-independent system $(\mathcal{Q},\mathcal{E},\mathcal{D})$  in Section~\ref{sec:ris}; cf.~\cite[Section~4.2.4.1]{mielke.roubicek}.
Precisely, Proposition~\ref{prop:transfer} asserts that for every {stable} sequence $(t_n,q_n)_{n\in\N}$ in $[0,T]\times \mathcal{Q}$ (i.e.,  each $q_n$ satisfies the global stability condition at time $t_n$)  such that  $(t_n,q_n)\to (t,q)$ in $[0,T]\times \mathcal{Q}$, and  every $\widetilde{q}\in\mathcal{Q}$, there exists a sequence $(\widetilde{q}_n)_{n\in\N}$ in $\mathcal{Q}$ satisfying
 $\widetilde{q}_n\to \widetilde{q}$ in $\mathcal{Q}$ and, notably,
\begin{equation*}
	\limsup_{n\to \infty} \big(  \mathcal{E}(t_n,\widetilde{q}_n)+\mathcal{D}(q_n,\widetilde{q}_n)  \big)\leq \mathcal{E}(t,\widetilde{q})+\mathcal{D}(q,\widetilde{q}).
\end{equation*}
This statement is exactly \cite[Eq.~(2.1.38)]{mielke.roubicek}, which, by \cite[Proposition~2.1.15(ii)]{mielke.roubicek},  implies that the {stability set} is closed \cite[Eq.~(C2)]{mielke.roubicek}. 

The proof of Proposition~\ref{prop:transfer} ultimately rests on the following lemma
which will allow us to move cavitation points while preserving the determinant constraint.

\begin{lemma}[Incompressible diffeomorphisms]\label{lem:incompressible}
Let ${A}\subset \Omega$   with $M\coloneqq\H^0({A})<+\infty$ and set ${A}=\{a_1,\dots,a_M\}$.
Then, there exists  ${\varepsilon}={\varepsilon}(\Omega,{A})>0$ such that
for every $B\subset\Omega$ with $\H^0(B)=M$, $B=\{b_1,\dots,b_M\}$, and
\begin{equation}\label{2610152254}
	\max \{|a_i-b_i|:\:i=1,\dots,M\}\leq {\varepsilon},
\end{equation}
there is $\varphi\in \mathrm{Diff}^\infty(\closure{\Omega};\closure{\Omega})$ with 
 \begin{equation}
 	\label{eq:h}
  	\varphi(b_i)=a_i  \quad \text{for all}\ i=1,\dots,M, \qquad \det D \varphi=1 \ \text{in}\ \Omega .
  \end{equation}
 Moreover, 
 \begin{equation}
 	\label{eq:h-est}
 	\|\varphi-\mathrm{id}\|_{C^1(\closure{\Omega};\RN)}\leq C(N) \sqrt{\max \{|a_i-b_i|:\:i=1,\dots,M\}}.
 \end{equation}
\end{lemma}


The proof of Lemma~\ref{lem:incompressible} rests on results about the prescribed Jacobian equation. Specifically, we employ an existence result for the corresponding boundary-value problem, which  we revise in Appendix~\ref{sec:ab}.

\begin{proof}[Proof of Lemma~\ref{lem:incompressible}]
	If ${A}=B$, then the claim is trivial as we can take $\varphi=\id$. Therefore, suppose  ${A}\neq B$. 
	The proof is subdivided into four steps.
	
	\emph{Step 1: Preparation.}
	Let
	\begin{equation*}
		d \coloneqq \min \{\dist(a_i,\partial \Omega):\:i=1,\dots,M \} \wedge \min \{ |a_i-a_j|:\:i,j=1,\dots,M,\:i\neq j \},
	\end{equation*}
	where the last minimum is understood as $+\infty$ when $M=1$, and fix $\eta\in C^\infty(\R;[0,1])$ such that
	\begin{equation}
		\label{eq:eta}
		\eta(t)=1\  \text{for all}\ t\leq 1/3,  \qquad  \eta(t)=0 \ \text{for all}\ t\geq 2/3.
	\end{equation}
	Set 
	\begin{equation*}
		\varepsilon\coloneqq \left( \frac{1}{16} \wedge \frac{d}{4} \wedge  \frac{ 1}{2K}   \right)^2, 
	\end{equation*}
	where $K\coloneqq \max \{ \|\eta\|_{L^\infty(\R)}, \|D \eta\|_{L^\infty(\R; \RN)}, \|D^2 \eta\|_{L^\infty(\R;\RNN)} \}$. The set $B$ is given in such a way that \eqref{2610152254} holds. Let
	\begin{equation}
		\label{eq:sigma}
		\delta\coloneqq \max \Big\{ \sqrt{|a_i-b_i|}:\:i=1,\dots,M \Big\}\leq \sqrt\varepsilon.
	\end{equation}
It is then possible to see that 
	\begin{equation*}
		\closure{B}(b_i,\delta)\cap \closure{B}(b_j,\delta)=\emptyset, \qquad \closure{B}(b_i,\delta)\subset \Omega, \qquad a_i\in \closure{B}(b_i,\delta/4)\subset B(b_i,\delta).
	\end{equation*}
Specifically this follows since 
	for all $i,j=1,\dots,\mu$ with $i\neq j$, we have
	\begin{equation*}
		|b_i-b_j|\geq |a_i-a_j|-|a_i-b_i|-|a_j-b_j|\geq d-2\delta^2\geq d - \delta/2>2\delta,
	\end{equation*}
	\begin{equation*}
		\dist(b_i,\partial \Omega)\geq \dist(a_i,\partial \Omega)-|a_i-b_i|\geq d-\delta^2  >2\delta,
	\end{equation*}
	and
	\begin{equation*}
		|a_i-b_i|\leq \delta^2 \leq \delta/4.
	\end{equation*}
	thanks to \eqref{eq:sigma}.
	
	\emph{Step 2: Local construction without incompressibility.}  
In this step we define a  diffeomorphism mapping $b_i$ to $a_i$.   This  transformation   does not  comply with   incompressibility   and   will be refined in the subsequent step of the proof. 
	Fix $i=1,\dots,M$. As in \cite[Lemma~9.1]{moracorral}, we define 
	$\widetilde{\varphi}_i \colon \closure{B}(b_i,\delta)\to \closure{B}(b_i,\delta)$ by setting
	\begin{equation*}
		\widetilde{\varphi}_i(x)\coloneqq x +  \eta\left(\frac{|x-b_i|}{\delta}\right)(a_i-b_i).
	\end{equation*} 
	In view of \eqref{eq:eta}, we have (cf.~Figure \ref{fig:varphii}):
	\begin{equation}
		\label{eq:h1}
		\widetilde{\varphi}_i=\mathrm{id}+a_i-b_i \ \text{in}\ \closure{B}(b_i,\delta/3),  \qquad \widetilde{\varphi}_i=\id \ \text{in}\ \closure{B}(b_i,\delta)\setminus B(b_i,2\delta/3).  
	\end{equation}
	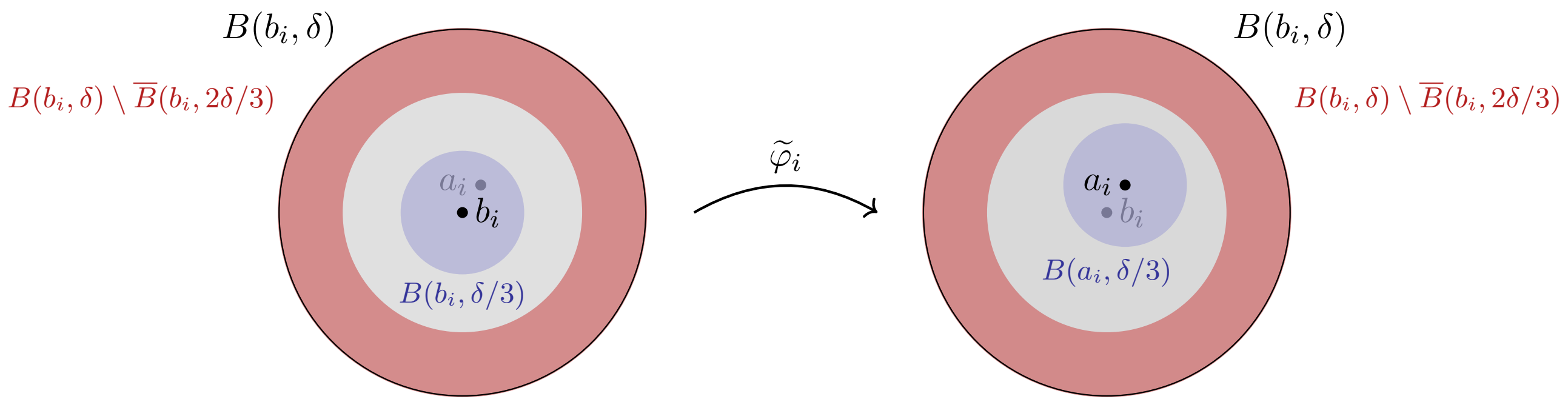
\begin{figure}
		\centering
		\begin{tikzpicture}
				\draw[fill=black!12,xshift=-100pt] (0,0) circle (2);
				\filldraw[black,xshift=-100pt] (.2,.3) circle (1.5pt) node[left] {$a_i$};
				\filldraw[bb!40,xshift=-100pt,opacity=.7,thick] (0,0) circle (.66);
				\filldraw[rr!60, even odd rule,xshift=-100pt,opacity=.8,thick]
				(0,0) circle (2)
				(0,0) circle (1.32);
				\draw[fill=black,xshift=-100pt] (0,0) circle (1.5pt) node[right] {$b_i$};
				\node[xshift=-100pt] at (-2,2) {$B(b_i,\delta)$};
				\node[xshift=-100pt,rr] at (-3.5,1.25) {\footnotesize $B(b_i,\delta)\setminus \closure{B}(b_i,2\delta/3)$};
				\node[xshift=-100pt,bb] at (0,-.9) {\footnotesize $B(b_i,{\delta}/{3})$};
				\draw[xshift=-100pt] (0,0) circle (2);

				\draw[->,thick,xshift=-85pt] (2,0) to[out=30,in=150] node[midway,above]  {$\widetilde{\varphi}_i$} (4,0);
				
				\draw[fill=black!15,xshift=100pt] (0,0) circle (2);
				\draw[fill=black,xshift=100pt] (0,0) circle (1.5pt) node[right ] {$b_i$};
				\filldraw[bb!40,xshift=100pt,opacity=.7,thick] (.2,.3) circle (.66);
				\filldraw[rr!60, even odd rule,xshift=100pt,opacity=.8,thick]
				(0,0) circle (2)
				(0,0) circle (1.32);
				\filldraw[black,xshift=100pt] (.2,.3) circle (1.5pt) node[left] {$a_i$};
					\node[xshift=100pt] at (2,2) {$B(b_i,\delta)$};
					\node[xshift=100pt,bb] at (0,-.65) {\footnotesize $B(a_i,{\delta}/{3})$};
				\node[xshift=100pt,rr] at (3.5,1.25) {\footnotesize $B(b_i,\delta)\setminus \closure{B}(b_i,2\delta/3)$};	
				\draw[xshift=100pt] (0,0) circle (2);
		\end{tikzpicture}
		\caption{The function $\widetilde{\varphi}_i$ in Step 2 of the proof of Lemma \ref{lem:incompressible}. }
		\label{fig:varphii}
	
	\end{figure}
	
	For all $x_1,x_2\in \closure{B}(b_i,\delta)$,
	\begin{equation*}
		\begin{split}
			|\widetilde{\varphi}_i(x_1)-\widetilde{\varphi}_i(x_2)-x_1+x_2|&\leq \left | \eta\left( \frac{|x_1-b_i|}{\delta}\right) - \eta\left( \frac{|x_2-b_i|}{\delta}\right) \right | \,|a_i-b_i|\\
			&\leq \|\eta'\|_{L^\infty(\R)}\delta|x_1-x_2|\leq K\delta|x_1-x_2|\leq \frac{1}{2}|x_1-x_2|
		\end{split}
	\end{equation*}
	thanks to \eqref{eq:sigma} and, in turn, $|\widetilde{\varphi}_i(x_1)-\widetilde{\varphi}_i(x_2)|\geq (1/2)|x_1-x_2|$.
	Thus, $\widetilde{\varphi}_i$ is injective and, actually, 
	$\widetilde{\varphi}_i\in \mathrm{Diff}^\infty(\closure{B}(b_i,\delta);\closure{B}(b_i,\delta))$. 
	We compute
	\begin{align*}
		D\widetilde{\varphi}_i(x)&=I+\frac{1}{\delta} \eta' \left(\frac{|x-b_i|}{\delta}\right)(a_i-b_i)\otimes \frac{x-b_i}{|x-b_i|},\\
		\det D \widetilde{\varphi}_i(x)&=1+\frac{1}{\delta}\eta' \left(\frac{|x-b_i|}{\delta}\right)(a_i-b_i)\cdot \frac{x-b_i}{|x-b_i|},\\
		D(\det D \widetilde{\varphi}_i)(x)&=\frac{1}{\delta^2} \eta''\left(\frac{|x-b_i|}{\delta}\right) \left(\frac{x-b_i}{|x-b_i|} \otimes \frac{x-b_i}{|x-b_i|} \right)(a_i-b_i)\\
		&\quad +\frac{1}{\delta |x-b_i|}\eta'\left(\frac{|x-b_i|}{\delta}\right)\left(I- \frac{x-b_i}{|x-b_i|} \otimes \frac{x-b_i}{|x-b_i|} \right)(a_i-b_i).
	\end{align*}
	Using \eqref{eq:sigma}, we estimate 
	\begin{align}
		\label{eq:h-tilde}
		\|\widetilde{\varphi}_i-\id\|_{C^0(\closure{B}(b_i,\delta);\RN)}&\leq K\delta^2,\\
		\label{eq:Dh-tilde}
		\|D\widetilde{\varphi}_i-I\|_{C^0(\closure{B}(b_i,\delta);\RNN)}&\leq K\delta\leq\frac{1}{2},\\
		\label{eq:det-Dh-tilde}
		\|\det D\widetilde{\varphi}_i-1\|_{C^0(\closure{B}(b_i,\delta))}&\leq  K\delta\leq\frac{1}{2}.
	\end{align}
	Meanwhile,  since $\eta$ is constant around zero, 
	\begin{align}
		\label{eq:Ddet-Dh-tilde}
		\|D(\det D \widetilde{\varphi}_i)\|_{C^0(\closure{B}(b_i,\delta);\RN)}\leq C(N).
	\end{align}
where $C(N)>0$ also depends on $K$; however we regard $\eta$ as a fixed function, so we neglect the dependence of the constants on $K$.

	\emph{Step 3: Local construction with incompressibility.} 
Now we have to modify the change of variable found in the previous step in such a way to guarantee incompressibility. We may argue separately for each ball. 
To work on a fixed domain, we introduce the  change of variable  $\zeta_i\colon \closure{B}(0,1)\to \closure{B}
	(b_i,\delta)$ given by $\zeta_i(p)\coloneqq b_i+\delta p$. The map
	$\widetilde{\psi}_i\coloneqq \zeta_i^{-1}\circ \widetilde{\varphi}_i\circ \zeta_i$ belongs to $   \mathrm{Diff}^\infty(\closure{B}(0,1);\closure{B}(0,1))$ and  satisfies
	\begin{equation}
		\label{eq:w}
		\widetilde{\psi}_i=\id + \frac{a_i-b_i}{\delta}  \ \text{in}\ \closure{B}(0,1/3), \qquad \widetilde{\psi}_i=\id \ \text{in}\ \closure{B}(0,1)\setminus {B}(0,2/3)
	\end{equation}
	in view of \eqref{eq:h1}.
	Moreover, from \eqref{eq:h-tilde}--\eqref{eq:Dh-tilde}, we obtain 
	\begin{align}
		\label{eq:wi-tilde}
		\|\widetilde{\psi}_i-\id\|_{C^1(\closure{B}(0,1);\RN)}&\leq \frac{1}{\delta}\|\widetilde{\varphi}_i-\id\|_{C^0(\closure{B}(b_i,\delta);\RN)}+\|D\widetilde{\varphi}_i-I\|_{C^0(\closure{B}(b_i,\delta);\RNN)}\leq 2K \delta.
	\end{align}
	Let $f_i\coloneqq \det D \widetilde{\psi}_i^{-1}\in C^\infty(\closure{B}(0,1))$. 
	Note that $f_i\geq 2/3$ in $B(0,1)$ 
	by the chain rule and  \eqref{eq:det-Dh-tilde}. Also, as $\widetilde{\psi}_i$ is a diffeomorphism, 
	\begin{equation*}
		\int_{B(0,1)}f_i\,\d  p =\int_{B(0,1)}\det D \widetilde{\psi}_i^{-1}\,\d  p = \leb(\widetilde{\psi}_i^{-1}(B(0,1)))=\leb(B(0,1))
	\end{equation*}
	thanks to the area formula.
	Using the chain rule and \eqref{eq:Dh-tilde}--\eqref{eq:det-Dh-tilde}, we find 
	\begin{equation}
		\label{eq:fi}
		\begin{split}
			\|f_i-1\|_{C^0(\closure{B}(0,1))}&=\left \| \frac{1-(\det D \widetilde{\psi}_i)\circ \widetilde{\psi}_i^{-1}}{(\det D \widetilde{\psi}_i)\circ \widetilde{\psi}_i^{-1}} \right \|_{C^0(\closure{B}(0,1))} \\
			&\quad \leq 2 \|\det D \widetilde{\varphi}-1\|_{C^0(\closure{B}(b_i,\delta);\RN)}\leq 2K \delta.
		\end{split}
	\end{equation}
	Then, 
	\begin{equation*}
		Df_i=\left( D \left[ \frac{1}{\det D \widetilde{\psi}_i} \right]\circ \widetilde{\psi}_i^{-1} \right) \left( (D\widetilde{\psi}_i)^{-1}\circ \widetilde{\psi}_i^{-1} \right),
	\end{equation*} 
	with
	\begin{equation*}
		\begin{split}
			D \left[ \frac{1}{\vphantom{\Big|}\det D \widetilde{\psi}_i} \right]&=-\frac{D \left[ \det D \widetilde{\psi}_i \right]}{\vphantom{\Big|}(\det D \widetilde{\psi}_i)^2}=-\frac{D\left[ (\det D \widetilde{\varphi}_i)\circ \zeta_i  \right]}{\vphantom{\Big|}(\det D \widetilde{\varphi}_i)^2\circ \zeta_i}=-\delta \left(\frac{ D\left[\det D \widetilde{\varphi}_i \right] }{\vphantom{\Big|} (\det D \widetilde{\varphi}_i)^2} \right) \circ \zeta_i
		\end{split}
	\end{equation*}
	and
	\begin{equation*}
		(D\widetilde{\psi}_i)^{-1}=(D\widetilde{\varphi}_i)^{-1}\circ \zeta_i=\sum_{k=0}^\infty (-1)^k (D\widetilde{\varphi}_i-I)^k\circ \zeta_i,
	\end{equation*}
	where we have used the Neumann series, whose convergence is guaranteed by \eqref{eq:Dh-tilde}. Thus,
	\begin{equation*}
		\label{eq:Dfi}
		\begin{split}
			\|Df_i\|_{C^0(\closure{B}(0,1);\RNN)}&\leq \delta \left \| \frac{ D\left[\det D \widetilde{\varphi}_i \right] }{\vphantom{\Big|} (\det D \widetilde{\varphi}_i)^2} \right \|_{C^0(\closure{B}(b_i,\delta);\RNN)} \|(D\widetilde{\varphi}_i)^{-1}\|_{C^0(\closure{B}(b_i,\delta);\RNN)}\\
			&\quad \leq 4\delta \left \|  D\left[\det D \widetilde{\varphi}_i \right]\right \|_{C^0(\closure{B}(b_i,\delta);\RNN)} \sum_{k=0}^\infty \|D\widetilde{\varphi}_i-I\|^k_{C^0(\closure{B}(b_i,\delta);\RNN)}
			\\
			&\quad \leq C(N)\,\delta,
		\end{split}
	\end{equation*} 
where the last inequality follows by  \eqref{eq:Dh-tilde} and \eqref{eq:Ddet-Dh-tilde}.
	This, together with \eqref{eq:fi}, yields
	\begin{equation}
		\label{eq:fi-est}
		\|f_i-1\|_{C^1(\closure{B}(0,1))}\leq C(N)\,\delta .
	\end{equation}
	Equation \eqref{eq:w} implies
	\begin{equation*}
		f_i=1 \ \text{in}\ \closure{B}\left(({a_i-b_i})/{\delta},1/3\right) \cup \left( \closure{B}(0,1)\setminus B(0,2/3) \right).
	\end{equation*}
	Note that $\closure{B}(0,1/4)\subset B((a_i-b_i)/\delta,1/3)$. Indeed, 
	\begin{equation*}
		\left |z-\frac{a_i-b_i}{\delta} \right |\leq |z|+\frac{|a_i-b_i|}{\delta}\leq \frac{1}{4}+\delta\leq \frac{1}{4}+\frac{1}{16}=\frac{5}{16}<\frac{1}{3}  
	\end{equation*}
	for each $z\in\RN$ with $|z|\leq 1/4$ thanks to \eqref{eq:sigma}.
	Therefore, we can apply Theorem \ref{thm:jac-det-eq} with
	\begin{equation*}
		\O= B(0,1),\quad  U= B(0,4/5)\setminus \closure{B}(0,1/5), \quad  V= B(0,3/4)\setminus \closure{B}(0,1/4), \quad  f= f_i
	\end{equation*}
	to find $v_i\in \mathrm{Diff}^\infty(\closure{B}(0,1);\closure{B}(0,1))$ satisfying
	\begin{equation}
		\label{eq:vi-eq}
		\begin{cases}
			\det D v_i=f_i \ \text{in}\ B(0,1), \\ \spt( v_i-\id)\subset U.
		\end{cases}
	\end{equation}
	Specifically, we deduce the estimate
	\begin{equation}
		\label{eq:vi}
		\|v_i-\id\|_{C^1(\closure{B}(0,1);\RN)}\leq C \|f_i-1\|_{C^1(\closure{B}(0,1))}\leq C(N)\,\delta
	\end{equation}
	owing to \eqref{eq:fi-est}.
	
	Define $\psi_i\coloneqq v_i\circ \widetilde{\psi}_i\in\mathrm{Diff}^\infty(\closure{B}(0,1);\closure{B}(0,1))$. From \eqref{eq:vi-eq}, we infer
	\begin{equation}
		\label{eq:wi-eq}
		\begin{cases}
			\det D \psi_i=1\quad \text{in}\ B(0,1), \\ \spt(\psi_i-\id)\subset B(0,4/5)
		\end{cases}
	\end{equation} 
	and, in particular,
	\begin{equation}
		\label{eq:wi0}
		\psi_i(0)=v_i(\widetilde{\psi}_i(0))=v_i\left(\frac{a_i-b_i}{\delta}\right)=\frac{a_i-b_i}{\delta}
	\end{equation}
	since
	\begin{equation*}
		\frac{|a_i-b_i|}{\delta}=\delta\leq \frac{1}{16}<\frac{1}{5}
	\end{equation*}
	by \eqref{eq:sigma}. Moreover, we claim that
	\begin{equation}
		\label{eq:wi-est}
		\|\psi_i-\id\|_{C^1(\closure{B}(0,1);\RN)}\leq C(N)\,\delta.
	\end{equation}
	Indeed,
	using \eqref{eq:wi-tilde} and \eqref{eq:vi}, we estimate
	\begin{equation*}
		\begin{split}
			\|\psi_i-\id\|_{C^1(\closure{B}(0,1);\RN)}&\leq \|\psi_i-\widetilde{\psi}_i\|_{C^1(\closure{B}(0,1);\RN)}+\|\widetilde{\psi}_i-\id\|_{C^1(\closure{B}(0,1);\RN)}\\
			&\leq\|v_i-\id\|_{C^0(\closure{B}(0,1);\RN)}\\
			&\quad +\|Dv_i-I \|_{C^0(\closure{B}(0,1);\RNN)}\|D\widetilde{\psi}_i\|_{C^0(\closure{B}(0,1);\RNN)}\\
			&\quad +\|\widetilde{\psi}_i-\id\|_{C^0(\closure{B}(0,1);\RN)} 
			+\|D\widetilde{\psi}_i-I\|_{C^0(\closure{B}(0,1);\RN)}\\
			&\leq C(N)\,\delta.
		\end{split}
	\end{equation*}
	Eventually, we change variables again and set $\varphi_i\coloneqq \zeta_i\circ \psi_i\circ \zeta_i^{-1}\in\mathrm{Diff}^\infty(\closure{B}(b_i,\delta);\closure{B}(b_i,\delta))$. The system \eqref{eq:wi-eq} gives
	\begin{equation}
		\label{eq:hi-eq}
		\begin{cases}
			\det D \varphi_i=1 \quad \text{in}\ B(b_i,\delta),\\ \spt(\varphi_i-\id)\subset B(b_i,4\delta/5),
		\end{cases}
	\end{equation}
	while \eqref{eq:wi0} immediately yields $\varphi_i(b_i)=a_i$. Also,
	\begin{equation}
		\label{eq:hi-est}
		\begin{split}
			\|\varphi_i-\id\|_{C^1(\closure{B}(b_i,\delta);\RN)}&\leq \delta\|\psi_i-\id\|_{C^0(\closure{B}(0,1);\RN)}+\|D\psi_i-I\|_{C^0(\closure{B}(0,1);\RNN)}\\
			&\leq C(N)\,\delta
		\end{split}
	\end{equation}
	thanks to \eqref{eq:wi-est}.
	
	\emph{Step 4: Global construction.} To conclude, it is sufficient to define $\varphi\colon \closure{\Omega}\to \closure{\Omega}$ by setting $\varphi\coloneqq \id$ in $\closure{\Omega}\setminus \bigcup_{i=1}^M B(b_i,\delta)$ and $\varphi\coloneqq \varphi_i$ in $B(b_i,\delta)$ for each $i=1,\dots,M$. Hence, from \eqref{eq:hi-eq}--\eqref{eq:hi-est}, we find that $\varphi\in\mathrm{Diff}^\infty(\closure{\Omega};\closure{\Omega})$ satisfies \eqref{eq:h}--\eqref{eq:h-est}.
\end{proof}

With Lemma~\ref{lem:incompressible} at hand, we can move to the proof of Proposition~\ref{prop:transfer}. 

\begin{proof}[Proof of Proposition~\ref{prop:transfer}]
	For convenience, we subdivide the proof into two steps.
	
	\emph{Step 1: Construction of the sequence.}
	For each $n\in\N$, define $\widetilde{A}_n\coloneqq A_n\cup ( \widetilde{A}\setminus A )$.
	We write the proof in the case where $A$ is a nonempty proper subset of $\widetilde{A}$, so
	we may write  $A\coloneqq \{a_1,\dots,a_\mu\}$ and $\widetilde{A}=\{a_1,\dots,a_\mu,b_1,\dots,b_{M-\mu} \}$ with $\mu\coloneqq \H^0(A)$ and $M\coloneqq \H^0(\widetilde{A})$. 
When $A=\widetilde A$ (so $\mu=M$), one may follow the very same proof, just dropping the points $b_j$. 
Up to subsequences, we write $A_n=\{a_{1,n}, \dots,a_{m,n}\}$ with $m\coloneqq \H^0(A_n)\geq \mu$ in such a way that
	\begin{enumerate}[\rm (i)]
		\item For each $i=1,\dots,m$, there exists a unique $j=1,\dots,\mu$ for which $a_{i,n}\to a_j$, as $n\to \infty$;
		\item For each $j=1,\dots,\mu$, there exists at least one $i=1,\dots,M$ for which $a_{i,n}\to a_j$, as $n\to \infty$;
	\end{enumerate}
see  \cite[Lemma~3.1]{bresciani.friedrich} for  details.
	In view of (ii), we select $\{i_1,\dots,i_\mu\}\subset \{1,\dots,m\}$ such that, $a_{i_j,n}\to a_j$, as $n\to \infty$, for all $j=1,\dots,\mu$.   
	Thus, for each $n\gg 1$, we can apply Lemma~\ref{lem:incompressible} with $\widetilde{A}$ in place of $A$ and  $B\coloneqq \{a_{i_1,n},\dots,a_{i_\mu,n},b_1,\dots,b_{M-\mu}  \}$ to find $\varphi_n\in\mathrm{Diff}^\infty(\closure{\Omega};\closure{\Omega})$ with 
	\begin{equation}
		\label{eq:phi-n-inc}
	\det D \varphi_n=1 \ \text{in}\ \Omega
	\end{equation}
	 such that
	\begin{equation}
		\label{eq:phi-n-points}
		\varphi_n(a_{i,n})=a_i \quad \text{for all}\ i=1,\dots,\mu, \qquad \varphi_n(b_k)=b_k \quad \text{for all}\ k=1,\dots,M-\mu
	\end{equation}
	and
	\begin{equation}
		\label{eq:phi-n-convergence}
		\|\varphi_n-\id\|_{C^1(\closure{\Omega};\RN)}\leq C(N)\sqrt{\max\{|a_{j}-a_{i_j,n}|:j=1,\dots,\mu  \}}.
	\end{equation}
	At this point, it is sufficient to set $\widetilde{z}_n\coloneqq \widetilde{z}\circ \varphi_n$ to have $(\widetilde{A}_n,\widetilde{z}_n)\in\widehat{\mathcal{C}}$. Indeed, $\det D \widetilde{z}_n=1$ almost everywhere in $\Omega$  and $C_{\widetilde{z}_n}=\varphi_n^{-1}(C_{\widetilde{z}})\subset \varphi_n^{-1}(\widetilde{A})\subset \widetilde{A}_n$ by \eqref{eq:phi-n-inc}--\eqref{eq:phi-n-points} and \cite[Proposition~5.2]{moracorral}.
	
	\emph{Step 2: Convergences.} Both the convergence of $(\widetilde{A}_n)_{n\in\N}$ in \eqref{eq:tc} and \eqref{eq:tI} are clear.
	Using the chain rule and \eqref{eq:phi-n-convergence}, we obtain the convergence of $(\widetilde{z}_n)_{n\in\N}$
	in \eqref{eq:tc}. Using again the chain rule with the change-of-variable formula and \eqref{eq:phi-n-inc}, we compute 
	\begin{equation*}
		\begin{split}
			\widehat{\mathcal{G}}(t_n,\widetilde{z}_n)&=\int_\Omega W\left( \Big( (Dw_{t_n})\circ \widetilde{z}\circ \varphi_n \Big) \Big( (D\widetilde{z})\circ \varphi_n \Big) D\varphi_n \right)\,\d x\\
			&=\int_\Omega W\left( \Big( (Dw_{t_n})\circ \widetilde{z} \Big)  D\widetilde{z} \Big( (D\varphi_n)\circ \varphi_n^{-1} \Big) \right)\,\d x.
		\end{split}
	\end{equation*} 
	Thus, the first limit in \eqref{eq:tGH} follows by dominated convergence in view of \eqref{eq:w-datum} and \eqref{eq:phi-n-convergence}. For the second limit, using \cite[Proposition~5.1]{moracorral} and change-of-variable formula for surface integrals, we write
	\begin{equation*}
		\begin{split}
			\widehat{\mathcal{H}}(t_n,\widetilde{z}_n)&=\sum_{a\in C_{\widetilde{z}_n}} \kappa(a) \per \big( w_{t_n}(E(\widetilde{z}_n,a)) \big)=\sum_{b\in C_{\widetilde{z}}} \kappa(\varphi_n^{-1}(b)) \per \big( w_{t_n}(E(\widetilde{z},b)) \big)\\
			&=\sum_{b\in C_{\widetilde{z}}} \kappa(\varphi_n^{-1}(b)) \int_{\partial^* E(\widetilde{z},b)} |(\cof D^{\partial^*E(\widetilde{z},b)}w_{t_n})\nu_{E(\widetilde{z},b)}|\,\d\haus.
		\end{split}
	\end{equation*}
	Letting $n\to \infty$, with the aid of \ref{it:kl-continuous} and \eqref{eq:phi-n-convergence}, the limits in \eqref{eq:tGH} follow once again by dominated convergence,
noting that each $t\mapsto w_t\restr{\partial^*E(\widetilde{z},b)}$ belongs to $W^{1,\infty}(0,T;W^{1,\infty}(\partial^*E(\widetilde{z},b);\Xi))$.
\end{proof}

\subsection{Characterization  of flaw points}\label{sec:flaw}
To complete the proof of Theorem~\ref{thm:auxiliary}, we are left to establish the characterization of flaw points in terms of cavitation points of deformations along the evolution.    

\begin{proposition}[Characterization of flaw points]
Let $(z,A)\colon [0,T]\to \widehat{\mathcal{C}}$ be a function satisfying conditions {\rm (i)--(iv)} in Theorem~\ref{thm:auxiliary}. Then, 
$A_t=A^0\cup \bigcup_{s\in [0,t]} C_{z_s}$ for all $t\in [0,T]$.
\end{proposition}

The proof of the proposition is quite elementary and works by contradiction.

\begin{proof}
Let $t\in (0,T]$ and set $B_t\coloneqq A^0 \cup \bigcup_{s\in [0,t]} C_{z_s}$. Condition (ii) in Theorem~\ref{thm:auxiliary} and the fact that $(z_t,A_t)\in \widehat{\mathcal{C}}$ yield the inclusion $B_t\subset A_t$. We show that $A_t\subset B_t$. By contradiction, let $a\in A_t\setminus B_t$. In particular,
\begin{equation}
	\label{eq:a}
	a\notin C_{z_s} \quad \text{for all}\ s\in [0,t].
\end{equation}
Let $t_0\coloneqq \inf \{s\leq t:\:a\in A_s   \}\leq t$. Suppose that $a\in A_{t_0}$. For a sequence $(s_n)_{n\in\N}$ in $(0,T)$ such that $s_n\to t_0^-$, we have $A_{s_n}\subset A_{t_0}\setminus \{a\}$ for all $n\in\N$. Owing to \eqref{eq:a}, $(z_{t_0},A_{t_0}\setminus \{a\})$ is an admissible competitor for the condition (iii) Theorem~\ref{thm:auxiliary} at each time $s_n$ and, in turn, 
\begin{equation*}
	\widehat{\mathcal{F}}(s_n,z_{s_n},A_{s_n})\leq \widehat{\mathcal{F}}(s_n,z_{t_0},A_{t_0}\setminus \{a\}) \quad \text{for all}\ n\in\N.
\end{equation*}
Given the continuity of $t\mapsto \widehat{\mathcal{F}}(t,z_t,A_t)$ claimed in condition (iv) of Theorem~\ref{thm:auxiliary} and \eqref{eq:w-datum}, letting $n\to \infty$ in the previous inequality yields
\begin{equation*}
	\widehat{\mathcal{F}}(t_0,z_{t_0},A_{t_0})\leq 	\widehat{\mathcal{F}}(t_0,z_{t_0},A_{t_0}\setminus \{a \})=	\widehat{\mathcal{F}}(t_0,z_{t_0},A_{t_0})-\lambda(a),
\end{equation*}
 which provides a contradiction since $\lambda(a)>0$ by \ref{it:kl-bdd}.

 Conversely, suppose that $a\notin A_{t_0}$. In this case,  we necessarily have $t_0<t$ since $a\in A_t$. Let $(s_n)_{n\in\N}$ be a sequence in $(0,T)$ with $s_n\to t_0^+$. Then,
 $s_n< t$ for $n\gg 1$.  By condition (ii) in Theorem \ref{thm:auxiliary}, we have $A_{t_0}\subset A_{s_n}\setminus \{a\}$. Also, \eqref{eq:a}  ensures that $(z_{s_n},A_{s_n}\setminus \{a\})$ is an admissible competitor for (iii) in Theorem~\ref{thm:auxiliary} at time $t_0$. Hence, 
 \begin{equation}
 	\label{eq:f}
 		\widehat{\mathcal{F}}(t_0,z_{t_0},A_{t_0})\leq \widehat{\mathcal{F}}(t_0,z_{s_n},A_{s_n}\setminus \{a\}).
 \end{equation} 
 Using Proposition~\ref{prop:td-energy}, the right-hand side is bounded as
 \begin{equation*}
 	\begin{split}
 		\widehat{\mathcal{F}}(t_0,z_{s_n},A_{s_n}\setminus \{a\})&\leq \widehat{\mathcal{F}}(s_n,z_{s_n},A_{s_n}\setminus \{a\})+ C(w)|s_n-t_0|\\
 		&= \widehat{\mathcal{F}}(s_n,z_{s_n},A_{s_n})-\lambda(a)+ C(w)|s_n-t_0|.
 	\end{split}
 \end{equation*}
 Hence, letting $n\to \infty$ in \eqref{eq:f}, once again by continuity of $t\mapsto \widehat{\mathcal{F}}(t,z_t,A_t)$  we obtain
 \begin{equation*}
 	\widehat{\mathcal{F}}(t_0,z_{t_0},A_{t_0})\leq 		\widehat{\mathcal{F}}(t_0,z_{t_0},A_{t_0})-\lambda(a),
 \end{equation*}
 which, as before, provides a contradiction. Therefore, $A_t\subset B_t$.
\end{proof}

\appendix

\section{First variation of the perimeter for bi-Lipschitz vector fields}
\label{sec:aa}

\setcounter{section}{0}
\setcounter{theorem}{0}

\setcounter{equation}{0}

\renewcommand{\thetheorem}{A.\arabic{theorem}}

\renewcommand{\theequation}{A.\arabic{equation}}

In this appendix, we revisit the derivation of the first-variation formula for the perimeter. Whereas classical treatments assume smooth deformation fields (see, e.g., \cite[Section~7.3]{ambrosio.fusco.pallara} or \cite[Chapter 7]{maggi}), here we consider fields that are merely bi-Lipschitz, in order to meet our assumption on Dirichlet data.
We demonstrate the following result, employed in Proposition \ref{prop:td-energy}.

\begin{theorem}[First variation of the perimeter]\label{thm:per}
	Let $T>0$, $\Xi\subset \RN$   the closure of a bounded Lipschitz domain, and $w\colon [0,T]\to W^{1,\infty}(\Xi;\Xi)$ satisfying \eqref{eq:w-datum}. 
	Then, for each $E\subset \Xi$ with $\per(E)<+\infty$, the function $P_E\colon [0,T]\to [0,+\infty)$ defined as $P_E(t)\coloneqq \per(w_t(E))$ belongs to $W^{1,\infty}(0,T)$. Moreover, there exists $S_w\subset (0,T)$ with $\mathscr{L}^1( S_w)=0$ for which the following holds:
	\begin{enumerate}[\rm (i)]
		\item For all $t\in (0,T)\setminus S_w$, there exists
		\begin{equation}
			\label{eq:PE-td}
			{P}_E'(t)= \int_{w_t(\partial^*E)} \div{}^{w_t(\partial^* E)} \left(\dot{w}_t\circ w_t^{-1} \right)\,\d \haus
		\end{equation}
		with
		\begin{equation}
			\label{eq:PE-est}
			|{P}_E'(t)|\leq C(N,w)  P_E(t).
		\end{equation}
		\item  There exists $\delta=\delta(w)>0$ and, for  each $t\in (0,T)\setminus S_w$, a function $\omega_t\colon (-\delta,\delta)\to [0,+\infty)$ with $\omega_t(t)\to 0$, as $h\to 0$, such that
		\begin{equation}
			\label{eq:PE-fo}
			\left | \frac{P_E(t+h)-P_E(t)}{h}-{P}_E'(t) \right |\leq C(N,w)\per(E)\omega_t(h) \quad \text{for all}\ h\in (-\delta,\delta).
		\end{equation}
	\end{enumerate}
\end{theorem}

The next result will be instrumental for the proof of Theorem~\ref{thm:per}.

\begin{lemma}[Differential estimates]\label{lem:Phi}
	Let $\nu\in \SN$ and define $\Phi_\nu\colon \RNN \to [0,+\infty)$  by setting $\Phi_\nu(F)\coloneqq |(\cof F)\nu|$. Also, let $F\in\RNN$ and $\ell>0$ be such that
	\begin{equation}
		\label{eq:Xlb}
		|F\tau|\geq \ell|\tau| \quad \text{for all $\tau \in \RN$ with $\tau\cdot \nu=0$}, \qquad F\nu=0. 
	\end{equation}
	Then, $\Phi_\nu$ is smooth at $F$. Moreover, for all $X,Y\in \RNN$, 
	\begin{align}
		\label{eq:dphi}
		\d \Phi_\nu(F)[X]&=\frac{(\cof F)\nu}{|(\cof F)\nu|}\cdot \big (\d(\cof)(F)[X]\big)\nu,\\ \label{eq:dphi2}
		\begin{split}
			\d^2\Phi_\nu(F)[X,Y]&=\frac{1}{|(\cof F)\nu|}\left(I-\frac{(\cof F)\nu}{|(\cof F)\nu|}\otimes \frac{(\cof F)\nu}{|(\cof F)\nu|}\right)\d(\cof)(F)[Y]\nu \cdot \d(\cof)(F)[X]\nu\\
			&\quad+ \frac{(\cof F)\nu}{|(\cof F)\nu|}\cdot \d^2(\cof)(F)[X,Y].
		\end{split}
	\end{align}	
	Additionally, for each $L>0$  such that $|F|\leq L$, we have
	\begin{equation}
		\label{eq:dPhi}
		|\d\Phi_\nu(F)[X]|\leq C(N,\ell,L)|X|, \qquad |\d^2\Phi_\nu(F)[X,Y]|\leq C(N,\ell,L)|X|\,|Y|.
	\end{equation}
\end{lemma}
\begin{proof}
First, we show that \eqref{eq:Xlb} entails 
\begin{equation}
	\label{eq:Philb}
	\Phi_\nu(F)\geq \ell^{N-1}.
\end{equation}
By assumption, $F$ has rank $N-1$. Let $\widetilde{\nu}\in\S^{N-1}$ be orthogonal to the range of $F$ and define $G\coloneqq F+\ell \widetilde{\nu}\otimes \nu$. For all $\xi\in\RN$, we have
\begin{equation*}
	|G\xi|^2=|F(I-\nu\otimes \nu)\xi|^2+\ell^2(\xi\cdot \nu)^2\geq \ell^2 \left( |(I-\nu\otimes \nu)\xi|^2+|\xi\cdot \nu|^2  \right)=\ell^2|\xi|^2.
\end{equation*}
thanks to \eqref{eq:Xlb} and, hence, $G$ is invertible with
\begin{equation}
	\label{eq:G1}
	|\det G|\geq \ell^N.	
\end{equation}
By the rank-one determinant identity and \eqref{eq:Xlb} once again,
\begin{equation*}
	\det G=\det F +\ell \widetilde{\nu}\cdot (\cof F)\nu=\ell \widetilde{\nu}\cdot (\cof F)\nu, 
\end{equation*}
which immediately gives
\begin{equation}
	\label{eq:G2}
	|\det G|\leq \ell |(\cof F)\nu|.
\end{equation}
The combination of \eqref{eq:G1}--\eqref{eq:G2} yields \eqref{eq:Philb}.

At this point, the smoothness of $\Phi_\nu$ at $F$ is clear and a direct computation yields \eqref{eq:dphi}--\eqref{eq:dphi2}. The entries of $\cof F$ are homogeneous polynomials of degree $N-1$ in the entries of $F$. In turn, for all $X,Y\in\RNN$,  there holds (at least for $N\geq 3$ as the cofactor is linear when $N=2$)
\begin{equation*}
	|\d(\cof)(F)[X]|\leq C(N)|F|^{N-2}|X|, \qquad |\d^2(\cof)(F)[X,Y]|\leq C(N) |F|^{N-3}|X|\,|Y|.
\end{equation*}
 Using these estimates and \eqref{eq:Philb} within \eqref{eq:dphi}--\eqref{eq:dphi2}, we easily deduce   
 \eqref{eq:dPhi}.
\end{proof}

\begin{proof}[Proof of Theorem \ref{thm:per}]
The proof is subdivided into four steps.  

\emph{Step 1: Preliminaries.} Fix $E\subset \Xi$ with $\per(E)<+\infty$. For brevity, we set $\Sigma\coloneqq \partial^*E$. Then, $\Sigma\subset \partial E\subset \Xi$ since  the latter set   is compact. The restriction on $\Sigma$ defines a nonexpansive linear operator from $\Lip_{\rm b}(\Xi;\RN)$ to $\Lip_{\rm b}(\Sigma;\RN)$ and, in turn, the function  $w\restr{\Sigma}\colon t\mapsto w_t\restr{\Sigma}$ belongs to $W^{1,\infty}(0,T;W^{1,\infty}(\Sigma;\Xi))$. Moreover, if $S_w\subset (0,T)$ with $\mathscr{L}^1(S_w)=0$ is such that $w$ is $W^{1,\infty}(\Xi;\RN)$-differentiable at all $t\in (0,T)\setminus S_w$, then $w\restr{\Sigma}$ is $W^{1,\infty}(\Sigma;\RN)$-differentiable at $t$ with derivative  
$\dot{w}_t\restr{\Sigma}$. 

By the Rademacher Theorem, each $w_t$ is $\haus$-almost everywhere tangentially differentiable with respect to $\Sigma$. More precisely, by arguing as in \cite[Remark~1.10]{lazzaroni}, we show that there exists  $H_w\subset \Sigma$ with $\haus(H_w)=0$ such that for all $t\in [0,T]$ and all $\xi\in \Sigma\setminus H_w$, $w_t$ is tangentially differentiable with respect to $\Sigma$ at $\xi$.

\emph{Step 2: Boundedness and Lipschitz continuity.} To simplify the notation, we let  $F_t\coloneqq D^{\Sigma}w_t$ and $\dot{F}_t\coloneqq D^{\Sigma}\dot{w}_t$.  Note that
\begin{equation}
	\label{eq:Ftan1}
	|F_t(\xi)\tau|\geq \ell |\tau|, \qquad F_t(\xi)\nu_E(\xi)=0 \qquad  \text{for all}\ t\in [0,T], \ \xi\in \Sigma\setminus H_w,\ \text{and}\ \tau\in \mathrm{T}_\xi\Sigma,
\end{equation}
and
\begin{equation}	
	\label{eq:Ftan2}
	\|F_t\|_{L^\infty(\Sigma;\RNN)}\leq L, \qquad \|F_t-F_s\|_{L^\infty(\Sigma;\RNN)}\leq L|t-s|\qquad \text{for all}\ s,t\in [0,T],
\end{equation}
where $\ell\coloneqq \|w^{-1}\|^{-1}_{W^{1,\infty}(0,T;W^{1,\infty}(\Xi;\Xi))}$ and $L\coloneqq \|w\|_{W^{1,\infty}(0,T;W^{1,\infty}(\Xi;\Xi))}$.
Now, define $\Phi_{\nu_E}$ as in Lemma~\ref{lem:Phi} and observe that the assumptions of the lemma are verified in view of \eqref{eq:Ftan1}--\eqref{eq:Ftan2}. 
By standard properties of sets of finite perimeter and the area formula,
\begin{equation}
	\label{eq:PE-cof}
	P_E(t)=\haus(\partial^*w_t(E))=\haus(w_t(\Sigma))=\int_{\Sigma} \Phi_{\nu_E}(F_t)\,\d\haus,
\end{equation}
where $\Phi_{\nu_E}(F_t)$ is the tangential Jacobian of $w_t$ with respect to $\Sigma$. Thus, $P_E$ is bounded since, by \eqref{eq:Ftan2},
\begin{equation*}
	|P_E(t)|\leq \int_\Sigma |\cof F_t|\,\d\haus \leq \haus(\Sigma) L^{N-1} \quad \text{for all}\ t\in [0,T].
\end{equation*}
Let $s,t\in [0,T]$ with $|t-s|\leq \ell/(2L)$. Fix  $\xi \in \Sigma\setminus H_w$. By  the Mean Value Theorem, 
\begin{equation*}
	\begin{split}
		|\Phi_{\nu_E(\xi)}(F_t(\xi))-\Phi_{\nu_E(\xi)}(F_s(\xi))|&\leq \int_0^1 \big|\d \Phi_{\nu_E(\xi)}\big (F_s(\xi)+\vartheta (F_t(\xi)-F_s(\xi))\big)[F_t(\xi)-F_s(\xi)]\big|\,\d\vartheta\\
		&\leq C(N,\ell,L) |F_t(\xi)-F_s(\xi)|\leq C(N,\ell,L)|t-s|
	\end{split}
\end{equation*}
thanks to \eqref{eq:dPhi} and \eqref{eq:Ftan2}. Here, to apply Lemma~\ref{lem:Phi}, we noted that  
for every $\vartheta\in[0,1]$ 
\begin{equation*}
 \|F_s+\vartheta (F_t-F_s)\|_{L^\infty(\Sigma;\RNN)}\leq L + \frac\ell2
\end{equation*}
and 
\begin{equation*}
	|(F_s(\xi)+\vartheta (F_t(\xi)-F_s(\xi))\tau|\geq (\ell/2)|\tau|\quad  \text{for all}\  \tau \in \mathrm{T}_\xi\Sigma
\end{equation*}
 because of \eqref{eq:Ftan1}--\eqref{eq:Ftan2} and the bound on $|t-s|$. By partitioning  $[0,T]$ in intervals of length at most $\ell/(2L)$ and using telescopic sums, we thus find
 \begin{equation}
 	\label{eq:Phi-lip}
 	|\Phi_{\nu_E(\xi)}(F_t(\xi))-\Phi_{\nu_E(\xi)}(F_s(\xi))|\leq C(N,\ell,L)|t-s|\quad \text{for all $s,t\in [0,T]$.}
 \end{equation}  
Therefore
\begin{equation*}
	|P_E(t)-P_E(s)|\leq C(N,\ell,L) \haus(\Sigma)|t-s| \quad \text{for all}\ s,t\in [0,T],
\end{equation*}
so that $P_E\in \Lip_{\rm b}([0,T])$.

\emph{Step 3: Proof of (i).}
The integrand on the right-hand side of \eqref{eq:PE-cof} can  can be differentiated with respect to time  at all $t\in (0,T)\setminus S_w$ and $\xi \in \Sigma\setminus H_w$. Recalling \eqref{eq:Phi-lip}, with the aid of  the dominated convergence theorem we compute  
\begin{equation}
	\label{eq:dot-PE}
	{P}_E'(t)=\int_\Sigma \d \Phi_{\nu_E}(F_t)[\dot{F}_t]\,\d\haus.
\end{equation}
For each $t\in [0,T]$, the set $w_t(E)\subset \Xi$ has finite perimeter and $\haus(\partial^* w_t(E) \triangle w_t(\partial^* E) )=0$. In particular, 
\begin{equation}
	\label{eq:reducedboundaryimage}
	\nu_{w_t(E)}\circ w_t=\frac{(\cof F_t)\nu_E}{|(\cof F_t)\nu_E|}\quad \text{$\haus$-almost everywhere in $\Sigma\setminus H_w$.}
\end{equation}
Let $v_t\coloneqq \dot{w}_t\circ w_t^{-1}$ and $G_t\coloneqq D^{w_t(\Sigma)}v_t$. Then, $\dot{F}_t=(G_t\circ w_t) F_t$. 
In view of \eqref{eq:dphi} and \eqref{eq:reducedboundaryimage},  we find
\begin{equation*}
	\d \Phi_{\nu_E}(F_t)[\dot{F}_t]=\nu_{w_t(E)}\circ w_t\cdot \d(\cof)(F_t)[\dot{F}_t]\nu_E\quad \text{$\haus$-almost everywhere in $\Sigma\setminus H_w$}.
\end{equation*}
Writing $F_t+\vartheta\dot{F}_t=(I+\vartheta (G_t\circ w_t))F_t$, at all points in $\Sigma \setminus H_w$ we compute
\begin{equation*}
	\begin{split}
		\d(\cof)(F_t)[\dot{F}_t]&=\frac{\d}{\d \vartheta}\bigg|_{\vartheta=0} \cof (F_t+\vartheta\dot{F}_t)= \left( \frac{\d}{\d \vartheta}\bigg|_{\vartheta=0} \cof(I+\vartheta G_t\circ w_t) \right)(\cof F_t)\\
		&=\d(\cof)(I)[G_t\circ w_t](\cof F_t)=\left( \big((\tr G_t)\circ w_t\big)I-G_t^\top\circ w_t \right)(\cof F_t),
	\end{split}
\end{equation*}
where we used  $\cof(XY)=(\cof X)(\cof Y)$ and  $\d(\cof)(I)[X]=(\tr X)I-X^\top$. Thus, at $\haus$-almost all points in $\Sigma\setminus H_w$ we have 
\begin{equation*}
	\begin{split}
		\d \Phi_{\nu_E}(F_t)[\dot{F}_t]&=\nu_{w_t(E)}\circ w_t\cdot \left(\big  ((\tr G_t)\circ w_t\big)I-G_t^\top\circ w_t\right)\nu_{w_t(E)}\circ w_t |(\cof F_t)\nu_E|\\
				&=(\tr G_t)\circ w_t|(\cof F_t)\nu_E|
	\end{split}
\end{equation*}
 given that $\nu_{w_t(E)}\cdot G_t^\top \nu_{w_t(E)}=G_t\nu_{w_t(E)}\cdot \nu_{w_t(E)}=0$ by the definition of tangential gradient. 
  Having in mind the definition of tangential divergence, this expression substituted into \eqref{eq:dot-PE}   yields \eqref{eq:PE-td} since
\begin{equation*}
	{P}_E'(t)=\int_\Sigma (\tr G_t)\circ w_t|(\cof F_t)\nu_E|\,\d\haus = \int_{w_t(\Sigma)} \tr G_t\,\d\haus 
\end{equation*}
thanks to the area formula. The estimate \eqref{eq:PE-est} follows immediately once we note that  
\begin{equation*}
	\sup_{t\in (0,T)} \|G_t\|_{L^\infty(\Sigma_t;\RNN)}\leq \ell L
\end{equation*}
since $G_t=(\dot{F}_t\circ w_t) D^{\Sigma}w_t^{-1}$.

\emph{Step 4: Proof of (ii).} Up to renaming $S_w$, we can  assume that each element of $(0,T)\setminus S_w$ is a Lebesgue point of $\dot{w}$.
Let $t\in (0,T)\setminus S_w$. For  $h\in (t,T-t)$,  we have
\begin{equation}
	\label{eq:AP}
	\frac{P_E(t+h)-P_E(t)}{h}-{P}_E'(t)=\int_{\Sigma} \left( \frac{\Phi_{\nu_E}(F_{t+h})-\Phi_{\nu_E}(F_t)}{h}-\d\Phi_{\nu_E}(F_t)[\dot{F}_t]  \right)\,\d\haus.
\end{equation}
By \eqref{eq:Ftan2}, the field $X_h\coloneqq h^{-1}(F_{t+h}-F_t)$ satisfies  
\begin{equation}
	\label{eq:X}
	\|X_h\|_{L^\infty(\Sigma;\RNN)}\le\frac{1}{h} \|w_{t+h}-w_t\|_{W^{1,\infty}(\Xi;\Xi)}\leq L.
\end{equation}
Fix $\xi\in \Sigma\setminus H_w$. By Taylor expansion
\begin{equation*}
	\begin{split}
		\Phi_{\nu_E(\xi)}&(F_{t+h}(\xi))-\Phi_{\nu_E(\xi)}(F_t(\xi))=\Phi_{\nu_E(\xi)}(F_t(\xi)+hX_h(\xi))-\Phi_{\nu_E(\xi)}(F_t(\xi))\\
		&=h\d\Phi_{\nu_E(\xi)}(F_t(\xi))[X_h(\xi)]+h^2\int_0^1 (1-\vartheta)\d^2\Phi_{\nu_E(\xi)}(F_t+(\xi)\vartheta h X_h(\xi))[X_h(\xi),X_h(\xi)]\,\d\vartheta
	\end{split}
\end{equation*}
and, in turn,
\begin{equation}
	\label{eq:ts}
		\begin{split}
			\frac{\Phi_{\nu_E(\xi)}(F_{t+h}(\xi))-\Phi_{\nu_E(\xi)}(F_t(\xi))}{h}&-\d\Phi_{\nu_E(\xi)}(F_t(\xi))[\dot{F}_t(\xi)]\\
			&\hspace{-2.75cm}= \d \Phi_{\nu_E(\xi)}(F_t(\xi))[X_h(\xi)-\dot{F}_t(\xi)]\\
			 &\hspace{-2.3cm}+h \int_0^1 (1-\vartheta)\d^2\Phi_{\nu_E(\xi)}(F_t(\xi)+\vartheta hX_h(\xi))[X_h(\xi),X_h(\xi)]\,\d\vartheta.  
		\end{split}
\end{equation}
Here, similarly to Step 1, we observed that 
\begin{equation*}
	\|F_t+\vartheta hX_h\|_{L^\infty(\Sigma;\RNN)}\leq L+\frac{\ell}{2}
\end{equation*} 
and 
\begin{equation*}
	|(F_t(\xi)+\vartheta h X_h(\xi))\tau|\geq \frac{\ell}{2}|\tau|\quad \text{for all}\  \tau\in \mathrm{T}_\xi\Sigma
\end{equation*} 
provided $|h|\leq \ell/(2L)$. 
Using \eqref{eq:dPhi} and the Fundamental Theorem of Calculus, we estimate
\begin{equation*}
	\begin{split}
		\|\d \Phi_{\nu_E}(F_t)[X_h-\dot{F}_t]\|_{L^\infty(\Sigma;\RNN)}&\leq C(N,L) \|X_h-\dot{F}_t\|_{L^\infty(\Sigma;\RNN)}\\
		&\leq C(N,L)\left  \|\frac{w_{t+h}-w_t}{h}-\dot{w}_t \right \|_{W^{1,\infty}(\Xi;\Xi)}\\
		&\leq C(N,L)  {\omega}_{t}(h), 
	\end{split}
\end{equation*}
where
\begin{equation*}
	{\omega}_{t}(h)\coloneqq \frac{1}{h}\int_{t\wedge (t+h)}^{t \vee (t+h)} \|\dot{w}_{\vartheta}-\dot{w}_t \|_{W^{1,\infty}(\Xi;\Xi)}\,\d \vartheta.
\end{equation*}
With the aid of \eqref{eq:dPhi} and \eqref{eq:X}, the value of the reminder is bounded as
 \begin{equation*}
 	|h| \int_0^1 |\d^2\Phi_{\nu_E(\xi)}(F_t(\xi)+\vartheta X_h(\xi))[X_h(\xi),X_h(\xi)]|\,\d\vartheta\leq C(N,\ell,L)|h|.
 \end{equation*}  
 Therefore, from \eqref{eq:AP} and \eqref{eq:ts} we obtain
 \begin{equation*}
 	\left | \frac{P_E(t+h)-P_E(t)}{h}-{P}_E'(t) \right |\leq C(N,\ell,L) \haus(\Sigma) \left( {\omega}_{t}(h)+|h|   \right)
 \end{equation*}
 for all $h\in (-\ell/(2L),\ell/(2L))$. Since $\omega_t(h)\to 0$, as $h\to 0$, by the Lebesgue Differentiation Theorem, this proves \eqref{eq:PE-fo} up to renaming $\omega_t$. 
\end{proof}

\stepcounter{section}
\section{Prescribed Jacobian equation}
\label{sec:ab}

\setcounter{section}{0}
\setcounter{theorem}{0}
\setcounter{equation}{0}

\renewcommand{\thetheorem}{B.\arabic{theorem}}

\renewcommand{\theequation}{B.\arabic{equation}}

In this appendix, we  revise the existence theory for the prescribed Jacobian equation according to our needs. For  a comprehensive treatment, we refer to \cite{csato.dacorogna.kneuss}. Our approach is based on Moser’s flow method \cite{moser}. While  unable at retrieving the optimal regularity of solutions~\cite{dacorogna.moser}, the method can provide information about their support. In addition, the approach is quite elementary and avoids Schauder estimates, in contrast with \cite{dacorogna.moser}. 

We establish the following result, which complements \cite[Theorem~10.11]{csato.dacorogna.kneuss} with a stability estimate. The result was instrumental for the proof of Lemma \ref{lem:incompressible}.

\begin{theorem}[Jacobian determinant equation]\label{thm:jac-det-eq}
Let $k\in\N$, $\O,U,V\subset \RN$ domains, and $f\in C^k(\closure{\O})$ with
\begin{equation*}
	\O \ \text{bounded}, \qquad V\subset \subset U\subset \subset \O, \qquad \spt(f-1)\subset V, \qquad  f>0 \quad \text{in}\ \O, \quad \int_\O f\,\d  \leb =\leb(\O)
\end{equation*}
Then, there exists  $v\in \mathrm{Diff}^k(\closure{\O}; \closure{\O})$ satisfying
\begin{equation}
	\label{eq:jac-eq}
	\begin{cases}
		\det D v=f \quad  \text{in}\ \O, \\\spt(v-\id)\subset U.
	\end{cases}   
\end{equation}
Moreover, setting
\begin{equation}
	\label{eq:mL}
m:= \max_{\closure\O}f >0 , \qquad L:= \|f\|_{C^k(\closure{\O})} \ge 1,
\end{equation} 
we have
\begin{equation}
	\label{eq:jac-est}
	\|v-\id\|_{C^k(\closure{\O};\RN)}\leq C(k,\O,U,V,m,L) \|f-1\|_{C^k(\closure{\O})}.
\end{equation}
\end{theorem}

For the reader's convenience, we review the main steps of the proof, adapting the arguments in~\cite{csato.dacorogna.kneuss}. 
We first address the equation $\div u=f-1$, obtained by formally linearizing  $\det D v=f$ at the identity. The next result refines \cite[Proposition~9.7]{csato.dacorogna.kneuss} according to our needs.

\begin{proposition}[Divergence equation]\label{prop:div-eq}
	Let $k\in\N$, $\O,U,V\subset \RN$ domains, and  $g\in C^k(\closure{\O})$ with 
	\begin{equation*}
		\O \ \text{bounded,} \qquad V\subset \subset U\subset \subset \O,\qquad \spt g\subset V, \qquad \int_\O g\,\d \leb =0.
	\end{equation*}
	Then, there exists $u\in C^k(\closure{\O};\RN)$ satisfying
	\begin{equation}
		\label{eq:div-eq}
		\begin{cases}
			\div u=g  \quad \text{in}\ \O, \\  \spt u \subset U.
		\end{cases} 
	\end{equation}
	Moreover, 
	\begin{equation}
		\label{eq:div-est}
		\|u\|_{C^k(\closure{\O};\RN)}\leq C(k,U,V)\|g\|_{C^k(\closure{\O})}.
	\end{equation}
\end{proposition}

To prove Proposition \ref{prop:div-eq}, we follow the strategy in \cite[Section 9.3.2]{csato.dacorogna.kneuss} which builds upon two lemmas. These  correspond to \cite[Lemma 9.8 \& Lemma 9.9]{csato.dacorogna.kneuss} and are proved via the  same arguments. Here, we include explicit norm estimates for the constructed objects in terms of the given data. 

\begin{lemma}[Divergence equation in a cube]\label{lem:div-eq-cube}
Let $k\in \N$, $Q\subset \RN$ a cube, and $g\in C_{\rm c}^k(Q)$ with 
\begin{equation*}
	\int_Q g\,\d \leb=0.
\end{equation*} 
Then, there exists $u\in C^k_{\rm c}(Q;\RN)$ satisfying
\begin{equation*}
	\div u=g\quad \text{in $Q$.}
\end{equation*}
Moreover,
\begin{equation*}
	\|u\|_{C^k(\closure{Q};\RN)}\leq C(k,N) \|g\|_{C^k(\closure{Q}).}
\end{equation*}
\end{lemma}

\begin{lemma}[Cut-off construction]\label{lem:div-eq-cutoff}
Let $k\in \N$, $U,V\subset \RN$ domains and $g\in C^k(\closure{U})$ with
\begin{equation*}
	V\subset \subset U, \qquad \spt g\subset V,\qquad \int_U g\,\d  \leb =0.
\end{equation*}
Then, for each finite family $\{U_1,\dots,U_\mu\}$ of open sets satisfying
\begin{equation*}
	V\subset \subset \bigcup_{i=1}^\mu U_i\subset \subset U, \qquad U_i\cap V\neq \emptyset\quad \text{for all $i=1,\dots,\mu$,}
\end{equation*}
there exist $g_1,\dots,g_\mu\in C^k(\closure{U})$ with
\begin{equation}
	\label{eq:cutoff1}
	\spt g_i\subset U_i\cap V, \qquad \int_U g_i\,\d \leb =0\quad \text{for all $i=1,\dots,\mu$}
\end{equation}
such that
\begin{equation}
	\label{eq:cutoff2}
	g=\sum_{i=1}^{\mu} g_i\quad \text{in $\closure{U}$.}
\end{equation}
In particular, 
\begin{equation}
	\label{eq:cutoff3}
	\|g_i\|_{C^k(\closure{U})}\leq C(k,V,U_i)\|g\|_{C^k(\closure{U})} \quad \text{for all $i=1,\dots,\mu$.}
\end{equation}
\end{lemma}

\begin{proof}[Proof of Proposition \ref{prop:div-eq}]
Let  $U_1,\dots,U_\mu\subset \subset U$ be cubes such that $\closure{V}\subset \bigcup_{i=1}^\mu U_i$ and $U_i\cap V\neq \emptyset$ for all $i=1,\dots,\mu$. By Lemma \ref{lem:div-eq-cutoff},  we find $g_1,\dots,g_\mu\in C^k(\closure{ U})$ for which \eqref{eq:cutoff1}--\eqref{eq:cutoff3} hold true. For each $i=1,\dots,\mu$, we apply Lemma \ref{lem:div-eq-cube} to find 
$u_i\in C^k_{\rm c}({U}_i;\RN)$ satisfying
\begin{equation}
	\label{eq:BB3}
	\div u_i=g_i \quad \text{in}\ U_i
\end{equation}
and
\begin{equation}
	\label{eq:BB4}
	\|u_i\|_{C^k(\closure{Q}_i;\RN)}\leq C(k,U_i) \|g_i\|_{C^k(\closure{U}_i)}.
\end{equation}
Setting $u\coloneqq \sum_{i=1}^\mu u_i$ after extending the functions $u_i$ to the whole $\O$ by zero, the field $u$ satisfies \eqref{eq:div-eq}--\eqref{eq:div-est} in view of \eqref{eq:cutoff1}--\eqref{eq:BB4}.
\end{proof}

Before presenting the proof of Theorem~\ref{thm:jac-det-eq}, we mention the following result which summarizes the main features of the flow method. The statement is adapted from \cite[Theorem 12.1]{csato.dacorogna.kneuss}.

\begin{lemma}[Flow method]\label{lem:flow}
Let $k\in \N$, $\mathrm{O}$ a bounded domain, $b\in C^0([0,1]; C^k (\closure{\rm O};\RN))$. 
Then, there exists a unique solution $\Theta \in C^0([0,1]; C^k( \closure{\rm O};\closure{\rm O}))$ to the system
\begin{equation}
	\label{eq:ode}
	\begin{cases}
		\frac{\d}{\d t}\Theta_t=b_t\circ \Theta_t & \text{in $(0,1]$,}\\
		\Theta_0=\id. & {}
	\end{cases}
\end{equation}
Moreover, setting
\begin{equation*}
	R\coloneqq \|b\|_{C^0([0,1];C^k(\closure{\rm O}))}, 
\end{equation*}
 exists $C=C(k,\mathrm{O},R)>0$ such that
\begin{equation*}
	\|\Theta_t-\id\|_{C^k(\closure{\rm O};\RN)}\leq C \int_0^t \|{b}_s\|_{C^k(\closure{\rm O};\RN)}\,\d s \qquad \text{for all $t\in [0,1]$.}
\end{equation*}
Additionally, $\Theta_t\in \mathrm{Diff}^k(\closure{\rm O};\closure{\rm O})$ and $\Theta_t=\id$ on $\partial \mathrm{O}$ for all $t\in [0,1]$. More precisely,
\begin{equation*}
	\bigcup_{t\in [0,1]} \spt (\Theta_t-\id)\subset K.
\end{equation*}
whenever  there exists a compact set $K\subset \mathrm{O}$  such that
\begin{equation*}
	\bigcup_{t\in [0,1]} \spt b_t\subset K.
\end{equation*}
\end{lemma}

\begin{proof}[Proof of Theorem~\ref{thm:jac-det-eq}]
Fix $m,L>0$ as in \eqref{eq:mL}. By applying  Proposition~\ref{prop:div-eq} with $g= f-1$, we find $u\in C^k(\closure{\O};\RN)$ satisfying
\begin{equation}
	\label{eq:div-eq-bis}
	\begin{cases}
		\div u=f-1 & \text{in $\mathrm{O}$,}\\
		\spt u\subset U
	\end{cases}
\end{equation}
and
\begin{equation}
	\label{eq:div-est-bis}
	\|u\|_{C^k(\closure{\rm O};\RN)}\leq C(k,U,V)\|f-1\|_{C^k(\closure{\rm O})}.
\end{equation}
Define  $b\in C^0([0,1];C^k(\closure{\rm O};\RN))$ by setting
\begin{equation*}
	b_t\coloneqq \frac{u}{t+(1-t)f.}
\end{equation*}
Recalling \eqref{eq:mL}, we have $t+(1-t)f\geq 1\wedge m$ for all $t\in [0,1]$ and, in turn,  $b$ is well defined. 
 
Let $\Theta \in C^1([0,1];C^k(\closure{\rm O};\closure{\rm O}))$  the unique solution to \eqref{eq:ode}. We claim that $v\coloneqq \Theta_1$ fulfills both \eqref{eq:jac-eq} and \eqref{eq:jac-est}. 
First, we note that
\begin{equation}
	\label{eq:XY1}
	\bigcup_{t\in [0,1]} \spt b_t\subset \spt u\subset U
\end{equation}
in view of \eqref{eq:div-eq-bis}. Using the  Fa\`a di Bruno chain rule \cite{consantine.savits}   and \eqref{eq:mL}, for all $t\in [0,1]$, we estimate
\begin{equation*}
	\left \| \frac{1}{t+(1-t)f} \right \|_{C^k(\closure{\rm O})}\leq C(N,k) \sum_{j=0}^{k} \frac{\| f \|_{C^k(\closure{\rm O})}^j}{m^{j+1}}\leq C(N,k,m,L), 
\end{equation*}
so that
\begin{equation}
	\label{eq:XY2}
	\begin{split}
			\left \| b_t \right \|_{C^k(\closure{\rm O};\RN)}&\leq C(k,N) \|u\|_{C^k(\closure{\rm O};\RN)}\,\left \| \frac{1}{t+(1-t)f} \right \|_{C^k(\closure{\rm O})}\\
			&\leq C(k,N,m,L)\|u\|_{C^k(\closure{\rm O};\RN)}\\
			&\leq C(k,N,m,L)\|f-1\|_{C^k(\closure{\rm O})}
	\end{split}
\end{equation}
thanks to the submultiplicativity of the norm in $C^k(\closure{\rm O})$ and to \eqref{eq:div-est-bis}. Having established \eqref{eq:XY1}--\eqref{eq:XY2}, the second condition in \eqref{eq:jac-eq} and the estimate \eqref{eq:jac-est} follows at once by  Lemma \ref{lem:flow}. 

We are left to show that the equation in \eqref{eq:jac-eq} is satisfied. For the reader's convenience, we repeat the argument in~\cite{dacorogna.moser}. 
Define $h\in C^1([0,1];C^k(\closure{\rm O}))$ as $h_t\coloneqq(t+(1-t)f\circ \Theta_t)\det D \Theta_t$. It is sufficient to show that $h$ is constant with respect to $t$, since in that case 
\begin{equation*}
	\det D v_1=h_1=h_0=f\circ \Theta_0\det D \Theta_0=f
\end{equation*}
 in view of \eqref{eq:ode}. Using the chain rule and \eqref{eq:ode}, we compute
 \begin{equation*}
 	\begin{split}
 		\frac{\d}{\d t} \big( t+(1-t) f\circ \Theta_t   \big)=1-f\circ \Theta_t + (1-t) (D f)\circ \Theta_t \cdot b_t
 	\end{split}
 \end{equation*}
 and
 \begin{equation*}
 	\begin{split}
 		\frac{\d}{\d t} \left( \det D \Theta_t \right)&=(\cof D \Theta_t):\frac{\d}{\d t} D \Theta_t
 		=(\cof D \Theta_t):D (b_t\circ \Theta_t)\\
 		&=(\cof D \Theta_t):\Big( \big( (D b_t)\circ \Theta_t \big) D \Theta_t  \Big)=(\div b_t)\circ \Theta_t \det D \Theta_t. 
 	\end{split}
 \end{equation*}
Hence
\begin{equation}
	\label{eq:xZ}
	\frac{\d}{\d t} h_t= \Big( 1-f\circ \Theta_t + (1-t) (D f)\circ \Theta_t \cdot b_t + \big(  t+(1-t) f\circ \Theta_t  \big)  (\div b_t)\circ \Theta_t  \Big)\det D \Theta_t. 
\end{equation} 
Recalling the definition of $b_t$, we find
\begin{equation*}
	\div b_t=\frac{\div u}{t+(1-t)f}-\frac{1-t}{\big(t+(1-t)f\big)^2} (D f)\cdot u=\frac{\div u-(1-t)(D f)\cdot b_t}{t+(1-t)f}.
\end{equation*}
Substituting this  expression into \eqref{eq:xZ} gives
\begin{equation*}
	\frac{\d}{\d t} h_t= (1-f +\div u  )\circ \Theta_t \det D \Theta_t=0, 
\end{equation*}
where the last equality follows from  \eqref{eq:div-eq}.
\end{proof}

\section*{Acknowledgements}
The authors are members of GNAMPA -- INdAM (Gruppo Nazionale per l'Analisi Matematica, la Probabilità e le loro Applicazioni -- Istituto Nazionale di Alta Matematica). 
Part of this work was conducted during visits of 
M.B.\ to the Department of Mathematics and Computer Science ``Ulisse Dini'' at the University of Florence
and of G.L.\ to the Institute of Analysis at Johannes Kepler University Linz: the support of both institutions is gratefully acknowledged. 

\subsection*{AI statement}
AI tools were employed by the authors solely for editorial assistance to improve  readability and  language  of the text. No generative AI was involved in the formulation of the problem or  the development of the mathematical proofs.

\vskip 30pt
{
	}

\end{document}